\documentclass{article}
\usepackage[utf8]{inputenc} 
\usepackage[T1]{fontenc}    
\usepackage{tikz}
\usepackage[style=numeric, backend=biber, sorting=none]{biblatex}
\usepackage{algorithm}
\usepackage{algpseudocode}
\usepackage{float}
\usepackage{amssymb}
\usepackage{amsmath}
\usepackage{hyperref}
\usepackage{mathtools}
\usepackage{graphicx}
\usepackage{subcaption}
\usepackage{siunitx}
\usepackage{booktabs}
\usepackage{tikz-cd}
\usepackage{stix}
\usepackage{bm}
\usepackage{amsbsy}
\usepackage{amsthm}
\usepackage{arxiv}
\usepackage{xspace}
\usepackage{xurl}
\newcommand{\fenics}{\texttt{FEniCs}\xspace}

\newtheorem{remark}{Remark}

\newcommand{\velocity}{\bm{v}}

\newcommand{\stdbase}{{\bm{e}}}

\newcommand\norm[1]{\left\lVert#1\right\rVert}
\DeclareMathOperator*{\argmax}{arg\,max}
\DeclareMathOperator*{\argmin}{arg\,min}

\newcommand{\inflowBoundary}{\Gamma_{-}}
\newcommand{\characteristicBoundary}{\Gamma_{0}}
\newcommand{\outflowBoundary}{\Gamma_{+}}

\newcommand{\sourceIntensity}{\bm{\lambda}}

\newcommand{\misfit}{\bm{y}}
\newcommand{\misfiti[1]}{{y_{#1}}}

\newcommand{\measurement}{\bm{d}}

\newcommand{\normal}{\bm{n}}

\newcommand{\obsO}{\mathcal{B}} 
\newcommand{\pto}{\mathcal{F}} 
\newcommand{\pts}{\mathcal{K}} 
 
\newcommand{\mta}{\mathcal{Q}} 
\newcommand{\adjoint}{q} 
\newcommand{\Nobservations}{N_{\text{obs}}} 

\newcommand{\ndof}{{n_{\text{dof}}}}

\newcommand{\R}{\mathbb{R}}

\newcommand{\parameterContinuous}{m}

\newcommand{\measTend}{T_{\text{obs}}}

\newcommand{\ansatzSources}{\mathcal{S}} 
 
\newcommand{\ansatzSourcesContinuous}{\mathcal{S}_{}}

\newcommand{\observationforeqx}{x_i^\text{obs}}
\newcommand{\observationforeqt}{t^{\text{obs}}_i}
\newcommand{\observationforequation}{(\observationforeqt,\observationforeqx)}
\newcommand{\x}{x}

\newcommand{\xbr}{{\bm{x}}}

\newcommand{\up}{{\vec{u}^{n+1}}}
\newcommand{\un}{{\vec{u}^{n}}}
\newcommand{\pp}{{\vec{\adjoint}^{n+1}}}
\newcommand{\pn}{{\vec{\adjoint}^{n}}}
\newcommand{\ansatzSpace}{\mathcal{V}_h}

\newcommand{\fracSolidus}[2]{#1/#2\,}

\renewcommand{\vec}{\bm}

\newcommand{\parameterfield}{m}

\newcommand{\hmax}{h_\mathrm{max}}
\newcommand{\lu}{l_\mathrm{u}=5\hmax}
\newcommand{\ld}{l_\mathrm{d}=15\hmax}

\newcommand{\bx}{\boldsymbol{x}}
\newcommand{\cs}{c^\mathrm{s}}
\newcommand{\ce}{c^\mathrm{e}}
\newcommand{\bcs}{\mathbf{c}^\mathrm{s}}
\newcommand{\bce}{\mathbf{c}^\mathrm{e}}
\newcommand{\np}{n_\mathrm{p}}
\newcommand{\bPi}{\boldsymbol{\Pi}}
\newcommand{\bPip}{\boldsymbol{\Pi}_\mathrm{p}}
\newcommand{\bPim}{\boldsymbol{\Pi}_\mathrm{m}}
\newcommand{\tcr}[1]{\textcolor{black}{#1}}

\title{Rapid Identification of Moving Contaminant Sources Through Physics-Based Modelling}

\date{} 					

\author{ \href{https://orcid.org/0009-0009-6325-3578}{\includegraphics[scale=0.06]{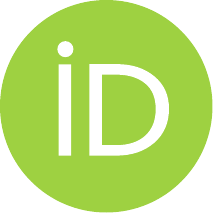}\hspace{1mm}Marco Mattuschka\textsuperscript{$1,$}}\thanks{Corresponding author e-mail: \texttt{marco.mattuschka@dlr.de}} ,\quad \href{https://orcid.org/0000-0001-8435-6466}{\includegraphics[scale=0.06]{orcid.pdf}\hspace{1mm}Jacopo Bonari\textsuperscript{$1$}},\quad
\href{https://orcid.org/0000-0002-2814-0027}{\includegraphics[scale=0.06]{orcid.pdf}\hspace{1mm}Max von Danwitz\textsuperscript{$1$}}, \quad
\href{https://orcid.org/0000-0002-8820-466X}{\includegraphics[scale=0.06]{orcid.pdf}\hspace{1mm}Alexander Popp\textsuperscript{$1,2$}}\\
\textsuperscript{$1$}German Aerospace Center (DLR), Institute for the Protection of Terrestrial Infrastructures,\\ 53757 Sankt Augustin, Germany\\
\textsuperscript{$2$}University of the Bundeswehr Munich, Institute for Mathematics and Computer-Based Simulation (IMCS),\\ 85577 Neubiberg, Germany
}

\renewcommand{\headeright}{}
\renewcommand{\undertitle}{}
\renewcommand{\shorttitle}{Algorithmic Identification of Moving Contaminant Sources}

\begin{document}

\maketitle

\begin{abstract}
In an act of sabotage or terrorism, hazardous material might be released deliberately into the atmosphere to threaten individuals, e.g., those operating critical infrastructure. Hazardous materials in such a scenario include toxic industrial chemicals (TICs), which are often invisible to the human eye, making it difficult to detect and respond to releases in a timely manner. This contribution considers the scenario of an airborne hazardous release requiring rapid and reliable assessment, with a chemical, biological, radiological, and nuclear (CBRN) sensor system providing scarce and local measurements. We present a novel algorithm that couples these data with an advection–diffusion model to detect, localize, and quantify a moving and time-varying contaminant source. Unlike many existing methods, the approach identifies sources with unknown occurrence time and trajectory by incorporating spatial sparsity as prior information. The feasibility of the approach is demonstrated in a two-dimensional computational domain. To further increase the technology readiness level, we additionally propose a calibration methodology for the required three-dimensional flow models based on wind tunnel experiments. Finally, a strategy for coupling the framework with real-time sensor data within a digital twin environment is outlined to enable predictive decision support in emergency scenarios.
\end{abstract}


\section{Introduction}

Airborne transport of hazardous substances poses a serious threat to communities and critical infrastructure. Releases of contaminants may occur accidentally, for example due to industrial leaks or spills, or intentionally in acts of sabotage or terrorism \cite{Boris.2002,Patnaik.2012,Danwitz.2024}. In emergency situations, decision-makers require reliable and timely information about the current state of contamination in order to initiate appropriate countermeasures. 

Airborne contaminant monitoring relies primarily on sensor-based detection systems capable of identifying CBRN-agents. Common approaches include electrochemical sensors~\cite{MADADELAHI2025117099} for selective toxic gas detection and mass spectrometry-based systems~\cite{WANG2025180041}, which provide high sensitivity and specificity for the analysis of complex chemical mixtures. Both technologies provide highly localized measurements and therefore correspond to point observations in large-scale mathematical models. In contrast, hyperspectral optical detection systems~\cite{spie:133423fc396d96d6105c143c6a891646b5384334,10.1117/12.692922} deliver spatially distributed information over comparatively large areas. Despite their technological diversity, all these systems share a fundamental limitation: measurements are available only locally in space and over a restricted time interval. At the same time, contaminants are transported by complex flow fields, often dominated by advection and turbulent mixing, which leads to highly nontrivial dispersion patterns. Consequently, it is in general not possible to infer the global distribution of contaminants over a large domain of interest directly from discrete sensor data without the support of a physics-based simulation model.

The objective of this work is therefore to extract actionable knowledge from spatio-temporally discrete measurements by coupling them with a flow-based transport model. In particular, we aim at the algorithmic identification of moving and time-varying contaminant sources and at providing reliable predictions of the resulting contaminant fields. 

Considering the current state of research in contaminant source detection, we propose a method that transfers existing inverse modeling approaches towards practical emergency-response scenarios. A major limitation of many transient identification methods, e.g., ~\cite{Villa.2021,MATTUSCHKA2026118854,Casas.2019, Leykekhman.2020, Biccari.2023, Monge.2020}, is the assumption that the release time is known a priori. In realistic scenarios, for instance in the case of an intentional attack, this information is typically unavailable. The present work therefore addresses the challenging problem of identifying and predicting a transient source with unknown activation time.

For example, \cite{WANG2025180041} present a promising measurement system based on mass spectrometry that has been used to monitor ship emissions and to verify compliance with exhaust gas regulations. In that application, the measured concentration data were correlated with known ship trajectories and atmospheric conditions. If such trajectory information is not available and one attempts instead to reproduce the sensor measurements by testing all possible source location candidates, the resulting inverse problem becomes severely ill-posed. 

To obtain stable and physically meaningful solutions, prior information must be incorporated into the reconstruction process. In the approach proposed here, this is achieved by assuming that the underlying sources are sparse in space. This structural prior is enforced by an appropriate regularization term, and leads to a well-posed optimization problem capable of identifying moving contaminant sources from limited sensor data.

At this point, however, we must acknowledge that a sparsity-promoting regularization in space and time does not constitute the most appropriate mathematical model for the problem under consideration. In the present setting, this choice induces a systematic bias: later activation times are implicitly favored, and reconstructed source locations tend to be shifted to the last possible time instance to cause the measurement signal. This work explicitly highlights this limitation. Nevertheless, the numerical results demonstrate that, for a sufficiently dense sensor configuration, the method still yields meaningful reconstructions. Even when the number of available sensors is significantly reduced, the approach provides a rough estimate of the source location together with a reliable prediction of the contaminant evolution. In emergency scenarios, such information already represents a substantial gain in actionable knowledge, despite the identified modeling shortcomings.

\section{Mathematical Modeling and Source Identification Algorithm}

\begin{figure}
\begin{subfigure}{0.54\textwidth}
\includegraphics[width=0.90\linewidth]{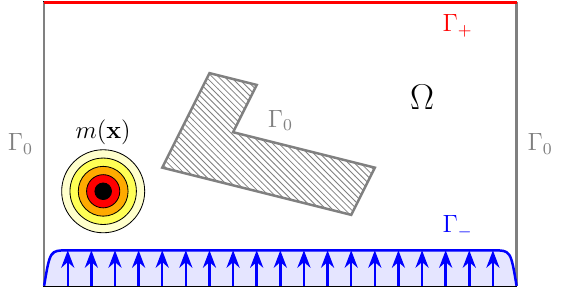}
\end{subfigure}
\begin{subfigure}{0.44\textwidth}
\includegraphics[width=0.80\linewidth]{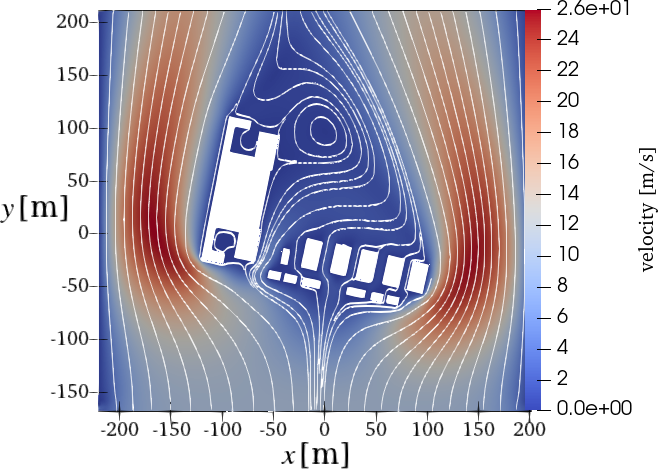}
\end{subfigure}
\caption{Computational domain $\Omega$ (left) with highlighted inflow ($\inflowBoundary$), outflow ($\outflowBoundary$), and characteristic ($\characteristicBoundary$) boundaries as well as the initial contaminant source $u(0,\cdot)=\parameterfield(\cdot)$ and the estimated wind vector field $\velocity$ (right) of benchmark scenario reproduced from~\cite{MattuschkaGoal}}
\label{fig:domain}
\end{figure}

\subsection{Numerical Modeling of Contaminant Sources and Dispersion}

\begin{figure}
\centering
\begin{subfigure}{0.48\textwidth}
\includegraphics[width=0.90\linewidth]{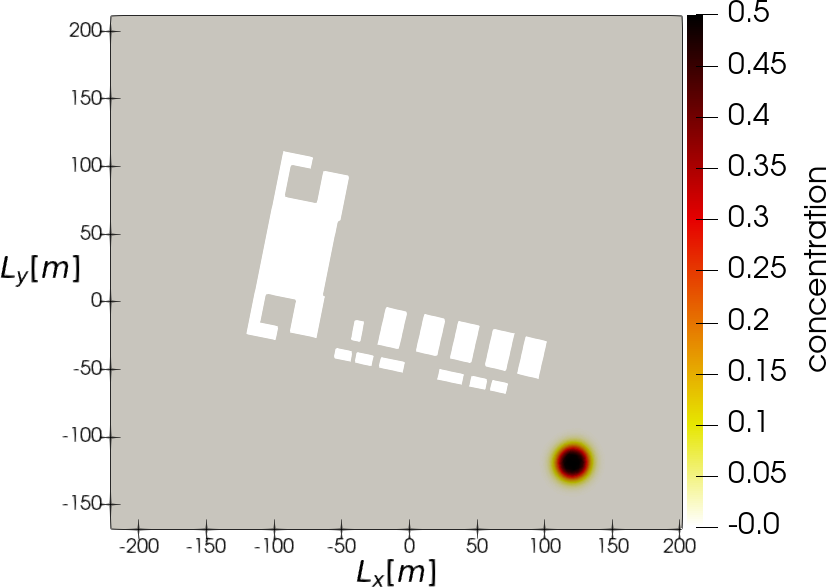}
\end{subfigure}
\begin{subfigure}{0.48\textwidth}
\includegraphics[width=0.90\linewidth]{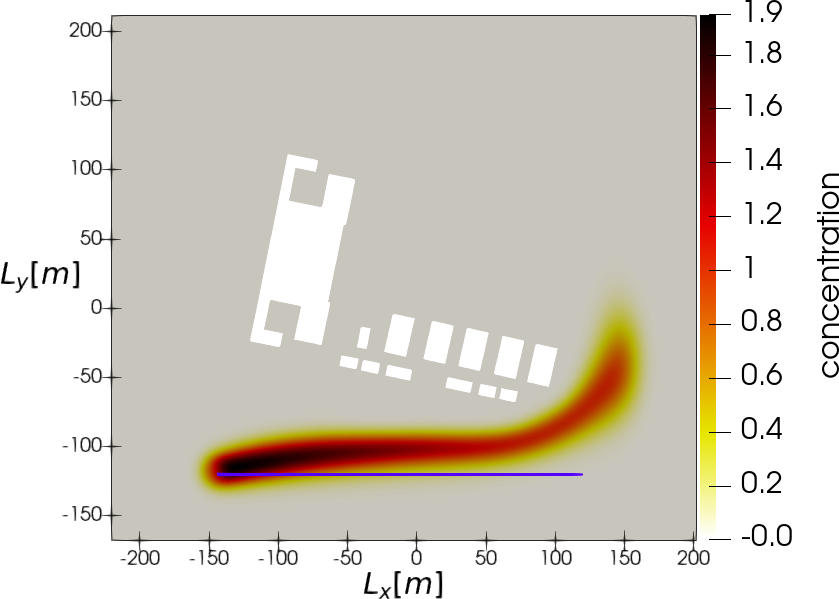}
\end{subfigure}
\begin{subfigure}{0.48\textwidth}
\includegraphics[width=0.90\linewidth]{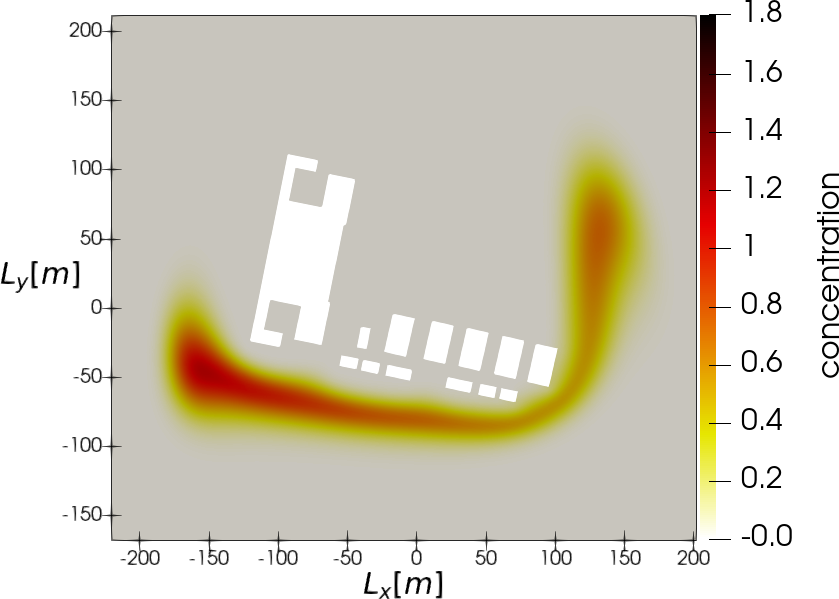}
\end{subfigure}
\begin{subfigure}{0.48\textwidth}
\includegraphics[width=0.90\linewidth]{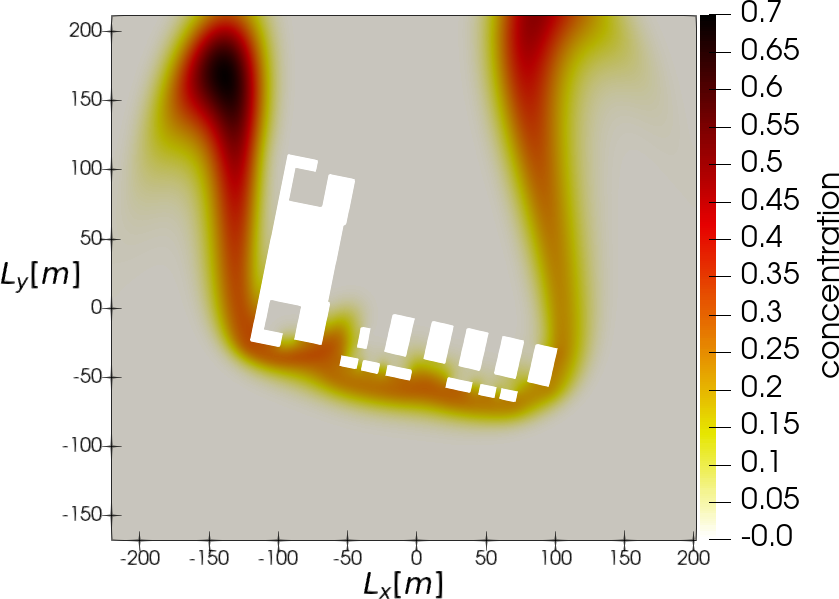}
\end{subfigure}
\caption{Numerical simulation of contamination dispersion at the beginning of the release ($t=\SI{0.5}{\s}$, top left), after completion of the release ($t=\SI{5.0}{\s}$) and the trajectory of the moving source (violet, top right), after the measurement period ($t=\SI{10.0}{\s}$, bottom left) and at the final simulation time ($t=\SI{20.0}{\s}$, bottom right).}
\label{fig:source}
\end{figure}

In the following, we investigate the dispersion of contaminants generated by a time-varying source. The source may move and vary in intensity over time. For the numerical simulations, we adopt an established benchmark scenario, cf.~\cite{Villa.2021,Danwitz.2024,MattuschkaGoal,MATTUSCHKA2026118854}. 

The underlying transport process is modeled with the advection-diffusion equation governed by the diffusion coefficient $\kappa$ and a wind vector field $\velocity$, which is assumed to be sufficiently smooth and divergence-free. The example wind field considered here is shown in \autoref{fig:domain}. Depending on the orientation of $\velocity$ relative to the outward-pointing boundary normal $\normal$, the boundary $\partial\Omega$ is decomposed into three disjoint subsets: the outflow boundary $\outflowBoundary \subset \partial\Omega$, where $\velocity \cdot \normal > 0$; the characteristic (tangential) boundary $\characteristicBoundary \subset \partial\Omega$, where $\velocity \cdot \normal = 0$; and the inflow boundary $\inflowBoundary \subset \partial\Omega$, where $\velocity \cdot \normal < 0$, following the convention in~\cite{Elman.2020}.

A mathematical description of the transport of a contaminant concentration $u$ in a bounded open domain $\Omega \subseteq \mathbb{R}^n$, $n \in \{2,3\}$, is given by the parameter-dependent forward problem
\begin{equation}\label{eq:forward_equation}
\begin{aligned}
(u_t - \kappa \Delta u + \velocity \cdot \nabla u)(t,x) 
&= \gamma_{\lambda}(t)\,\ansatzSources(\gamma_{x_s}(t),r,x) 
&& \text{in } (0,T)\times\Omega,\\
\kappa \nabla u \cdot \normal(t,x) 
&= 0 
&& \text{in } (0,T)\times (\outflowBoundary \cup \characteristicBoundary),\\
u(t,x) 
&= 0 
&& \text{in } (0,T)\times \inflowBoundary,\\
u(0,\cdot) 
&= 0 
&& \text{in } \Omega.
\end{aligned}
\tag{$\mathcal{P}_{\pts}$}
\end{equation}
To model realizations of a time-varying and moving source, we introduce a 
parameter curve 
\[
\gamma : [0,T] \rightarrow \R_{\geq 0} \times \Omega,
\qquad 
\gamma(t) = \bigl(\gamma_{\lambda}(t), \gamma_{x_s}(t)\bigr),
\]
where $\gamma_{\lambda}(t)$ denotes the source intensity and $\gamma_{x_s}(t)$ its spatial location at time $t$. In the present work, we restrict ourselves to contaminant sources whose spatial distribution is described by a smooth, radially symmetric shape 
function. Specifically, we consider the shape function
\begin{equation}
 \begin{aligned}
  \ansatzSources(\x_s,r,y)=\min \left\{0.5,\exp \left(-\text{ln}(\epsilon)\norm{y-\x_s}_2^2/r^2\right)\right\},
\end{aligned}\label{eq:ansatz}
\end{equation}
for the center $\x_s \in \Omega$, radius of the source $r>0$ and $\epsilon>0$ as a given threshold. Alternative realizations of the shape function $\ansatzSources$ can be 
found in~\cite{MATTUSCHKA2026118854}. 
An example of such a time-dependent source and the corresponding 
contaminant distribution obtained from 
\autoref{eq:forward_equation} is shown in \autoref{fig:source}.

The following describes how contaminant concentration measurements are obtained. To arrive at a formulation that is as close as possible to the actual application, we model sensor measurements as spatial averages rather than point measurements, in contrast to \cite{Villa.2021, MATTUSCHKA2026118854}. To this end, we introduce a smooth, radially symmetric basis function $\eta \colon \mathbb{R}^n \to \mathbb{R}$ such that $\eta(x)=1$ for $x \in B_{\sigma}(0)$ and $\operatorname{supp}(\eta) \subset B_1(0)$.

\begin{figure}
\begin{subfigure}{0.48\textwidth}
\includegraphics[width=0.90\linewidth]{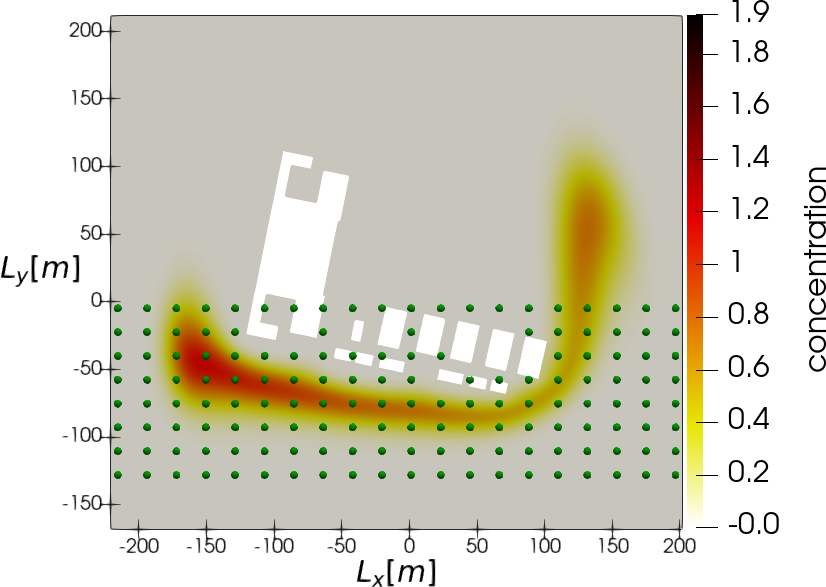}
\end{subfigure}
\begin{subfigure}{0.48\textwidth}
\includegraphics[width=0.90\linewidth]{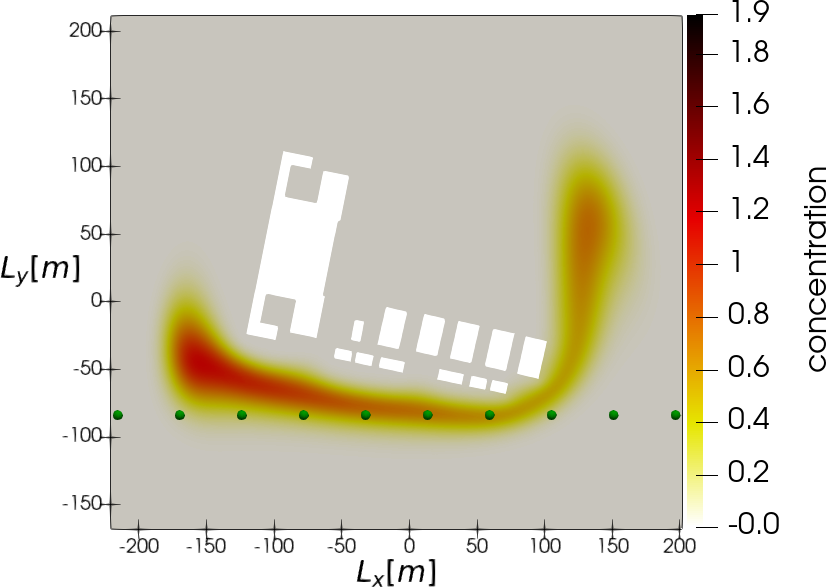}
\end{subfigure}
\caption{Numerical simulation of contaminant dispersion after the measurement period ($t=\SI{10}{\s}$) illustrated along with the dense sensor configuration (left) and the sparse sensor configuration (right). Sensor locations are marked by green dots.}
\label{fig:eom}
\end{figure}

For a solution $u$ of \autoref{eq:forward_equation}, the sensor measurement at a 
space-time location $(\observationforeqx,\observationforeqt) \in \Omega \times [0,\measTend]$, with $\measTend \leq T$, 
is defined as
\[
\misfiti[i] := \int_0^T \int_{\Omega} u(y,t)\, 
\eta(y-\observationforeqx, t-\observationforeqt)\, 
\mathrm{d}\Omega(y)\, \mathrm{d}t.
\]
This definition gives rise to a linear and bounded space-time \textit{observation operator}
\[
\obsO(u) := \sum_{i=1}^{\Nobservations} y_i\, \stdbase_i,
\quad \text{for a sequence of observation points } 
\observationforequation,\ i=1,\dots,\Nobservations.
\]
Consequently, we define the \textit{parameter-to-observable} operator by $\pts(\gamma)=u$ and 
$\pto(\gamma)=\obsO \circ \pts(\gamma)$. 
Given a \textit{misfit} vector $\misfit \in \mathbb{R}^{\Nobservations}$, for example
$
\misfit = \frac{1}{\sigma^2}\bigl(\obsO(u)-\measurement\bigr),
$
where the measurements $\measurement$ are perturbed by white noise, the associated 
\textit{misfit-to-adjoint} map is defined by $\mta(\misfit)=\adjoint$. 
Here, $\adjoint$ denotes the solution of the final value problem
\begin{equation}\label{eq:adjoint_equation}
\begin{aligned}
  -\adjoint_t - \kappa \Delta \adjoint - \operatorname{div}(\adjoint \velocity) 
  &= \sum_{i=1}^{\Nobservations} \misfiti[i] \, \eta_i
  && \text{in } (0,T)\times\Omega, \\
  (\velocity \adjoint + \kappa \nabla \adjoint) \cdot \normal 
  &= 0
  && \text{on } (0,T)\times(\outflowBoundary \cup \characteristicBoundary), \\
  \adjoint 
  &= 0
  && \text{on } (0,T)\times\inflowBoundary, \\
  \adjoint(T,\cdot) 
  &= 0
  && \text{in } \Omega,
\end{aligned}
\tag{$\mathcal{P}_{\mta}$}
\end{equation}
with the smooth right-hand side $\eta_i = \eta(y-\observationforeqx, t-\observationforeqt)$.

\begin{remark}
By replacing point-wise measurements with spatially averaged observations, \autoref{eq:adjoint_equation} admits a smooth solution. In earlier formulations, including \cite{Villa.2021,Danwitz.2024,MattuschkaGoal,MATTUSCHKA2026118854}, point evaluations were employed, which introduce Dirac distributions on the right-hand side of \autoref{eq:adjoint_equation}. When the unknown parameter corresponds to the initial condition, sufficient regularity can be recovered by starting the measurement process at a time $T_0 > 0$, exploiting the strong smoothing properties of the Laplacian. However, this imposes a restriction on the admissible observation window or modeled initial condition. The present formulation avoids this limitation by employing spatially averaged measurements, thereby ensuring regularity of the adjoint solution without requiring a delayed measurement start. 
\end{remark}

\subsection{A Sparse Inversion Framework for Source Detection}
Before formulating the inverse problem for identifying a transient contaminant source, 
we describe the discretization of the partial differential equations introduced in the 
first section. We employ a standard finite element discretization of the 
advection--diffusion problem, which yields the following discrete counterpart of 
\autoref{eq:forward_equation}:
\begin{equation}\label{eq:forward_equation_d}
\begin{aligned}
(M + \Delta t\, V + \Delta t\, \kappa K + \Delta t\, \tau S + \tau V^\top)\, \up
&= (M + \tau V^\top)\,(\un + \Delta t\, \vec{m}^{n}), \\
\vec{u}_{n=0} &= 0,
\end{aligned}
\end{equation}
where we use continuous Lagrange nodal basis functions defined by
\[
\ansatzSpace = \operatorname{span}\{\phi_1, \dots, \phi_{\ndof}\}
\]
associated with the nodes $\{p_1, \dots, p_{\ndof}\}$. 
The mass matrix, stiffness matrix, and skew-symmetric advection matrix are given by
\[
M_{ij} := \int_{\Omega} \phi_i(x)\,\phi_j(x)\,\mathrm{d}\Omega(x), \,
K_{ij} := \int_{\Omega} 
\langle \nabla \phi_i(x), \nabla \phi_j(x) \rangle\,\mathrm{d}\Omega(x),\,
V_{ij} := \int_{\Omega} 
\phi_i(x)\,\langle \nabla \phi_j(x), \velocity \rangle\,\mathrm{d}\Omega(x).
\]
An implicit Euler time discretization is performed using the approximation
\[
\vec{u}_t \approx \fracSolidus{(\up - \un)}{\Delta t}
\]
at time instances $(0, \Delta t, \dots, T = \Delta t\, N_T)$, yielding a solution in the discrete space--time space $\bigoplus_{i=0}^{N_T} \ansatzSpace$. The implementation of the discrete observation operator $\obsO^h: \bigoplus_{i=0}^{N_T} \ansatzSpace \rightarrow \mathbb{R}^{\Nobservations}$ works straightforwardly by projecting the function $\eta$ into the finite element space and evaluating the function by multiplying it with the mass matrix. In summary, this gives us the discrete counterpart ${\pto}^h:\bigoplus_{i=0}^{N_T} \ansatzSpace  \rightarrow \mathbb{R}^{\Nobservations} $ for a given right-hand side $\vec{\parameterContinuous}=(\vec{m}^0,\dots,\vec{m}^{N_T})$.

We employ the well-established SUPG stabilization 
technique \cite{Brooks.1982, Danwitz.2023}, which introduces the matrix
\[
S_{ij} := \int_{\Omega} 
\langle \nabla \phi_i(x), \velocity \rangle
\left( 
\langle \nabla \phi_j(x), \velocity \rangle 
- \kappa\, \Delta \phi_j(x) 
\right)
\,\mathrm{d}\Omega(x),
\]
together with the stabilization parameter
\[
\tau = \min\left( \fracSolidus{h_E^2}{2\kappa}, \fracSolidus{h_E}{\|\velocity\|} \right),
\]
where $h_E := \sup_{x,y \in E} |x-y|$ denotes the diameter of a finite element $E$.

In a very similar form, the discretization of the adjoint \autoref{eq:adjoint_equation} is given by
\begin{align}
\label{eq:dis_adj}
(M + \Delta \bar t V^T + \Delta \bar t \kappa K + \Delta \bar t \tau S^T + \tau V) {\pn} &= (M + \tau V) {\pp} + M \vec{\misfit}^{n+1}\,,
\end{align}
with initial \( \adjoint_{n=0} = 0 \). This variant leads to the discrete operator \( \mta^h : \mathbb{R}^{\Nobservations} \rightarrow \bigoplus_{i=0}^{n_T} \ansatzSpace \). For each time step $n \in \{0,\dots,N_T\}$, we denote by $\mta^h_n$ the projection onto the $n$-th component of this direct sum. The component $\mta^h_n$ corresponds to the solution of the discrete adjoint problem \autoref{eq:adjoint_equation} at time $t_n = n \Delta \bar t$.

To model the parameter $\parameterContinuous(t,x):=\gamma_{\lambda}(t)\,\ansatzSources(\gamma_{x_s}(t),r,x)$ in $\Omega$, we employ the 
standard $L^2$-projection of the source term from 
\autoref{eq:ansatz} onto the finite element space. In the discrete setting, 
this corresponds to the projection with respect to the mass matrix $M$, i.e.,
\[
\ansatzSources^h(\x_s,r,\cdot)
=
\argmin_{f \in \ansatzSpace}
\|f - \ansatzSources(\x_s,r,\cdot)\|_{M}^2.
\]
Its representation in coefficient form is denoted by the finite element 
vector $\vec{\ansatzSources}^h(\x_s,r)$. Given a parameter curve $\gamma : [0,T] \to \R_{\ge 0} \times \Omega,$ we define the discrete transient parameter $\vec{\parameterContinuous}$ at time levels $t_n = \Delta t\, n$ by
\begin{equation*}
\vec{\parameterContinuous}^{\,n}
= 
\gamma_{\lambda}(t_n)\,
\vec{\ansatzSources}^h\left(\gamma_{x_s}(t_n),\, r
\right),
\qquad
n \in \{1,\dots,N_T\}.
\end{equation*}
To stay consistent with the modeling of \cite{MATTUSCHKA2026118854}, we extend the framework by searching for a sparse representation of the source at each discrete point in time, i.e., sums of atoms represented as integrals of the shape function $\ansatzSources$. This means that we search for each time step $\Delta t\,n$, $N^n$ source locations $\xbr^{n} \in \bar{\Omega}^{N^n_{}}$ and intensities $\sourceIntensity^{n} \in \R^{N^n_{}}_{>0 }$.
In summary, this results in the following form on the right side of \autoref{eq:forward_equation_d} 
\begin{equation*}
 \vec{\parameterContinuous}^n=\sum^{N^n_{}}_{j=1}\,\lambda^{n}_j\,\ansatzSourcesContinuous(\x^{n}_j,\cdot), \, \text{where } (\xbr^{n},\sourceIntensity^{n}) \in \bar{\Omega}^{N^n_{}} \times \R^{N^n_{}}_{>0 } \text{ for each } n \in \{1,\dots,N_T\}
\end{equation*}
and the inverse problem now consists of minimizing the following functional
\begin{equation} \label{eq:objective}
 \min_{\parameterContinuous} \left\lbrack \fracSolidus{1}{(2\,\sigma^2)}\,\|{\pto}^h(\vec{\parameterContinuous}) - \measurement\|^2_{\R^{\Nobservations}} + \alpha\, \sum_{n=0}^{N_T}|{\sourceIntensity}^{n}|_{\ell_1} \right\rbrack  \text{ for all admissible sources } \vec{\parameterContinuous}=(\vec{m}^0,\dots,\vec{m}^{N^T}).
 \tag{$\mathcal{P}(\vec{m})$}
\end{equation}

\subsection{Algorithm for Source Tracking}
For the existence and uniqueness of a minimizer of 
\autoref{eq:finitealgo}, it is necessary to extend the space of admissible 
sources in \autoref{eq:objective} to the convex cone of positive Radon measures. 
For a detailed discussion of the modeling aspects and the corresponding proofs, 
we refer to~\cite{MATTUSCHKA2026118854,Pieper.2021}. 
As a consequence, \autoref{eq:objective} is minimized only over a finite set of candidates. More precisely, given a set of candidate locations 
$\vec{x}=\left(\xbr^{0},\dots,\xbr^{N_T}\right)$ and corresponding source 
intensities 
$\sourceIntensity=\left(\sourceIntensity^{0},\dots,\sourceIntensity^{N_T}\right)$, 
we define
\begin{equation}\label{eq:parametertosourceterm}
\vec{\parameterContinuous}[\vec{x},\sourceIntensity]
=
\left(
\sum^{N^0}_{j=1}\lambda^{0}_j\,\vec{\ansatzSourcesContinuous}(\x^{0}_j),
\dots,
\sum^{N^{N_T}}_{j=1}\lambda^{N_T}_j\,\vec{\ansatzSourcesContinuous}(\x^{N_T}_j)
\right).
\end{equation}
The finite-dimensional objective
\begin{equation} \label{eq:finitealgo}
 \min_{\sourceIntensity}
 \left[
 \frac{1}{2\sigma^2}
 \left\|
 \pto^h\bigl(\vec{\parameterContinuous}[\vec{x},\sourceIntensity]\bigr)
 - \measurement
 \right\|^2_{\R^{\Nobservations}}
 +
 \alpha \sum_{n=0}^{N_T}
 \left|
 \sourceIntensity^{n}
 \right|_{\ell_1}
 \right]
 \tag{$\mathcal{P}(\vec{x})$}
\end{equation}
is then minimized. This problem can be solved using established 
semi-smooth Newton methods \cite{Milzarek.2014}. The candidate locations are determined by extracting the maxima of the field
\begin{equation} \label{eq:dual}
{\varphi}^k_n
=
-
\int_\Omega
\ansatzSourcesContinuous(\cdot,z)\,
\mta({\misfit}_k)(\Delta t\,n,z)
\,\mathrm{d}\Omega(z)
\end{equation}
in each iteration. The entire algorithm is described in detail in \autoref{alg:pdap}. 

\begin{algorithm}
\caption{Primal-Dual-Active-Point-Strategy for Time-Varying Source Identification}\label{alg:pdap}
\begin{algorithmic}
\Require Nodal points $\left\{p_1, \dots, p_{\ndof}\right\}$, shape function $\vec{\ansatzSourcesContinuous}$, $\vec\parameterContinuous_0=(0,\dots,0)$, empty matrices and vectors ${\vec{x}_0=\left(\xbr^{0}_0,\dots,\xbr^{N^T}_0\right),\vec{\sourceIntensity_0}=\left(\sourceIntensity^{0}_0,\dots,\sourceIntensity^{N^T}_0\right)}$\\
\For{$k = 0,1,2...$}
\State~
\State 1.\,Given 
$\vec{\parameterContinuous}_k=m[\vec{x}_k,\vec{\lambda}_k]$ according to \autoref{eq:parametertosourceterm}, compute misfit$
  {\misfit}_k =  \fracSolidus{1}{\sigma^2} \left(\pto^h(\vec{\parameterContinuous}_k) - \measurement\right).$
\State~
\State 2.\,Compute convolution (cf. \autoref{eq:dual})
\begin{equation*}
 \vec{\varphi}^n_{k} =-
\left((M\vec{\ansatzSourcesContinuous}(p_0))^T \, \mta_0^h({\misfit}_k),\dots,(M\vec{\ansatzSourcesContinuous}(p_\ndof))^T \, \mta_{N_T}^h({\misfit}_k)\right)
\end{equation*}
\State 3.\,Determine index of the maximum
\begin{equation*}
    i^n_k\in \argmax \vec{\varphi}^n_{k} \text{ for each } n \in \{1,\dots,N_T\}.
\end{equation*}
\State 4.\,Set $\xbr_{k+1/2}=\xbr_{k}$. 
\State~
\State 5.\,\textbf{for }{$n = 0,\dots,N_T$}\\ 
\State~
\qquad \qquad Append $p_{i^n_k}$ to $\xbr^{n}_{k+1/2}$ \textbf{if} $ (\vec\varphi_{k}^n)_{i^n_k}> \alpha+ \text{tol}$\\
\State~
\,\textbf{if} no point has been added, \Return
\State \,\,\quad 
\State 6.\, Solve the finite minimization problem and update source intensities
\begin{equation*}
\left({\sourceIntensity}_{k+1/2}\right) \in \argmin \left(
\text{\hyperref[eq:finitealgo]{{$\mathcal{P}\left(\vec{x}_{k+1/2}\right)$}}}\right).
\end{equation*}
\State \,\,\quad
\State 7. Update iterates
\begin{equation*}
\vec{\parameterContinuous}_k=\vec{\parameterContinuous}_k[\vec{x}_{k+1/2},\sourceIntensity_{k+1/2}].
\end{equation*}
\State \,\,\quad
\State 8.\, Obtain $\xbr_{k+1},{\sourceIntensity}_{k+1}$ by pruning all entries for which ${\sourceIntensity}_{k+1}$ is negligibly small.
\State~
\EndFor
\end{algorithmic}
\end{algorithm}

\section{Numerical Simulation as Proof of Concept - Moving Contaminant Source on Campus}

\begin{figure}
\centering
\begin{subfigure}{0.48\textwidth}
\includegraphics[width=0.90\linewidth]{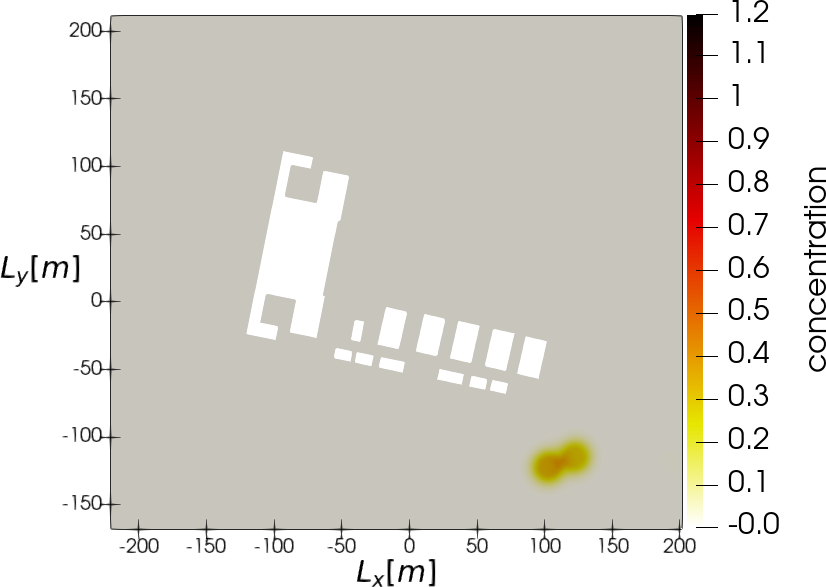}
\end{subfigure}
\begin{subfigure}{0.48\textwidth}
\includegraphics[width=0.90\linewidth]{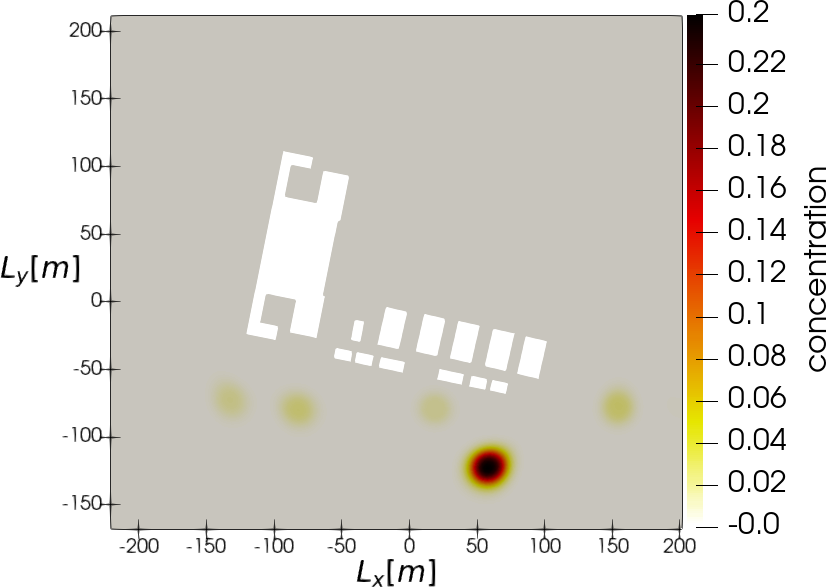}
\end{subfigure}
\caption{Reconstruction of the contaminant sources at the beginning of the release ($t=\SI{0.5}{\s}$) for the two considered sensor configurations (dense left, sparse right).}
\label{fig:predstart}
\end{figure}

\begin{figure}
\centering
\begin{subfigure}{0.48\textwidth}
\includegraphics[width=0.90\linewidth]{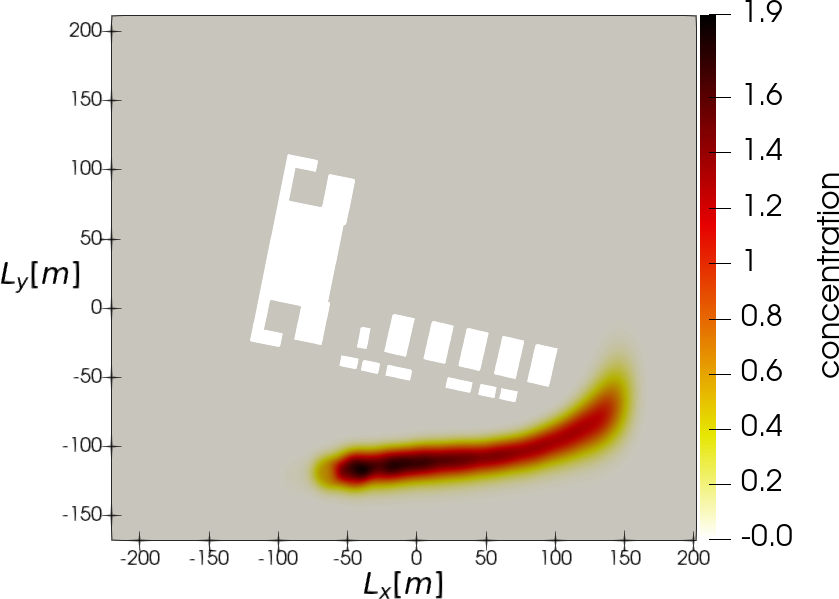}
\end{subfigure}
\begin{subfigure}{0.48\textwidth}
\includegraphics[width=0.90\linewidth]{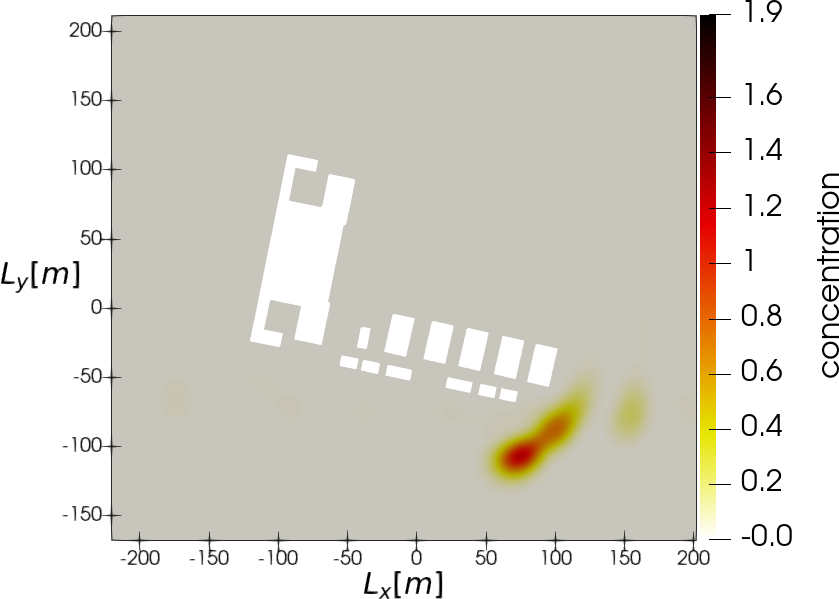}
\end{subfigure}
\caption{Reconstruction of the contaminant sources at $t=\SI{3.5}{\s}$ for the two considered sensor configurations (dense left, sparse right).}
\label{fig:pred04}
\end{figure}
\begin{figure}
\centering
\begin{subfigure}{0.48\textwidth}
\includegraphics[width=0.90\linewidth]{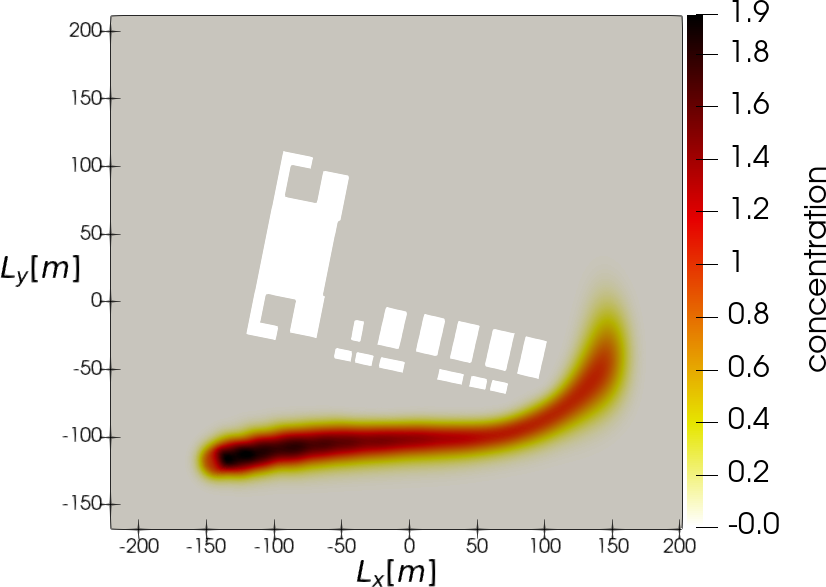}
\end{subfigure}
\begin{subfigure}{0.48\textwidth}
\includegraphics[width=0.90\linewidth]{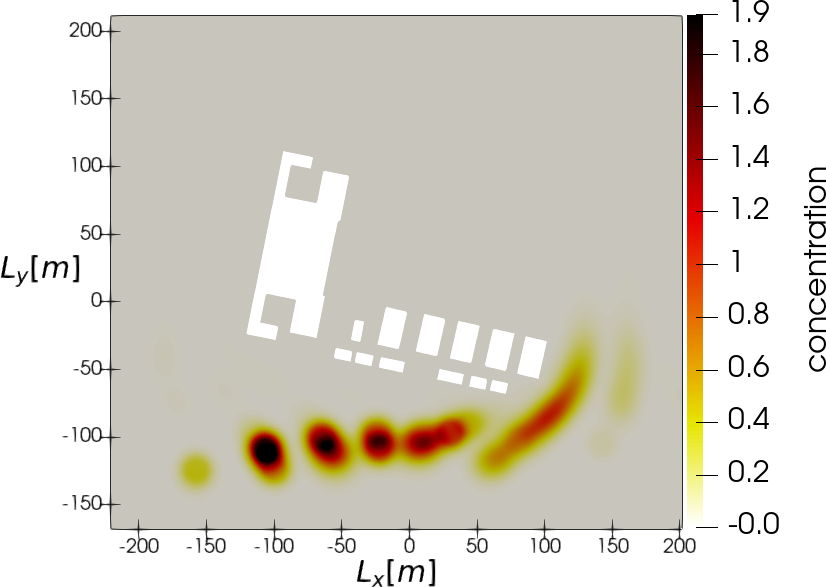}
\end{subfigure}
\caption{Reconstruction of the contaminant sources at the end of the release ($t=\SI{5.0}{\s}$) for the two considered sensor configurations (dense left, sparse right).}
\label{fig:eosource}
\end{figure}
\begin{figure}
\centering
\begin{subfigure}{0.48\textwidth}
\includegraphics[width=0.90\linewidth]{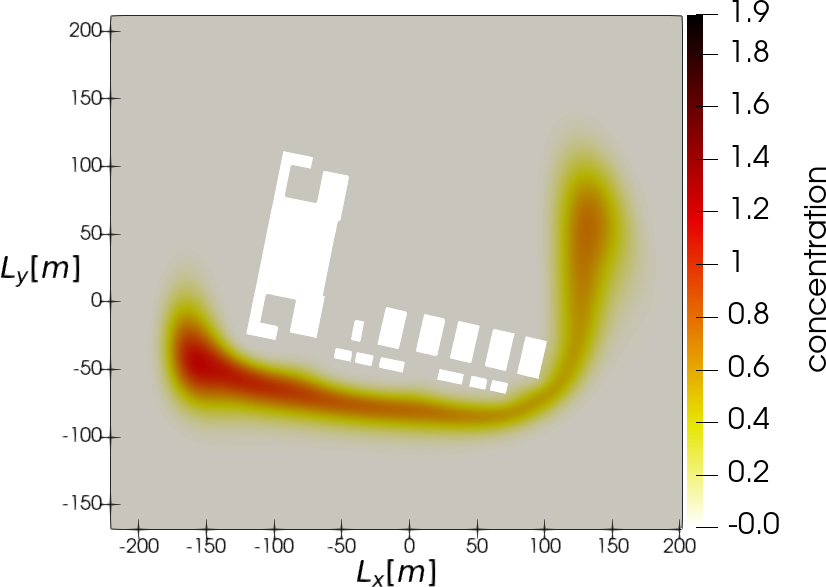}
\end{subfigure}
\begin{subfigure}{0.48\textwidth}
\includegraphics[width=0.90\linewidth]{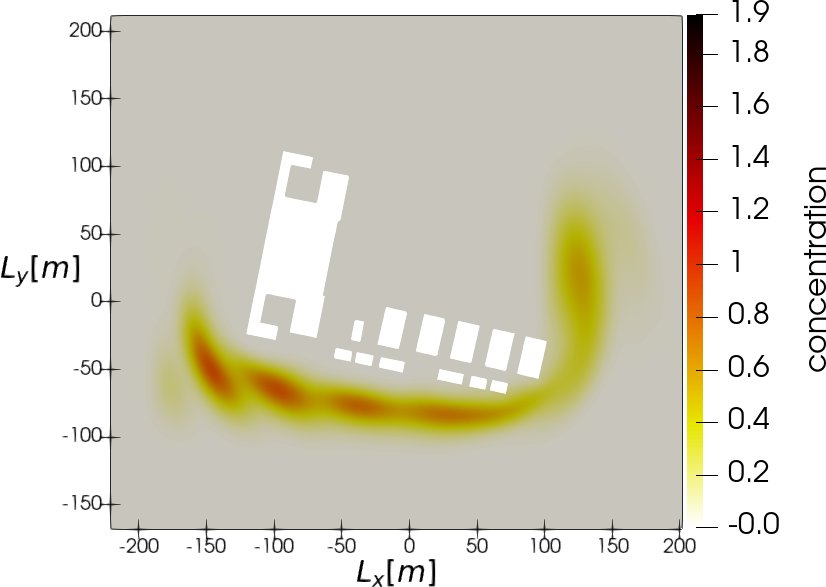}
\end{subfigure}
\caption{Reconstruction of the contaminant dispersion at the end of the measurement period ($t=\SI{10}{\s}$) for the two considered sensor configurations (dense left, sparse right).}
\label{fig:adjoint}
\end{figure}
\begin{figure}
\centering
\begin{subfigure}{0.48\textwidth}
\includegraphics[width=0.90\linewidth]{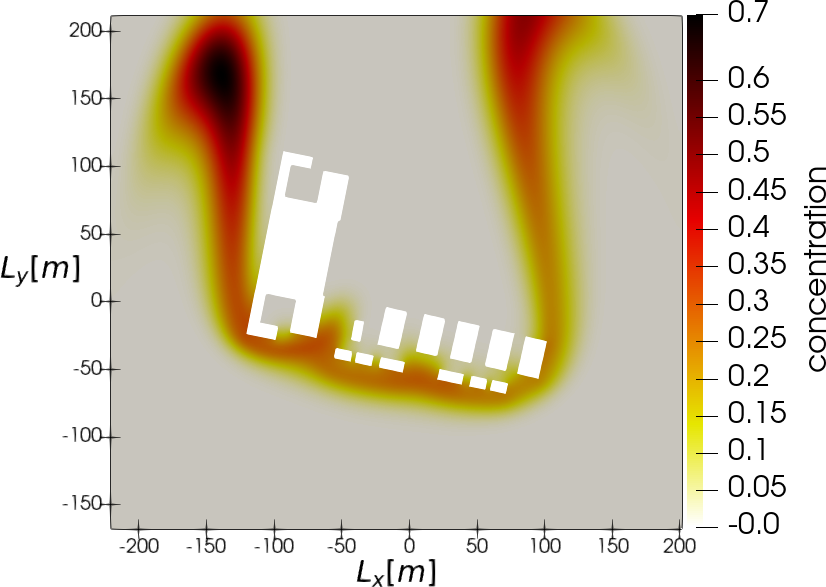}
\end{subfigure}
\begin{subfigure}{0.48\textwidth}
\includegraphics[width=0.90\linewidth]{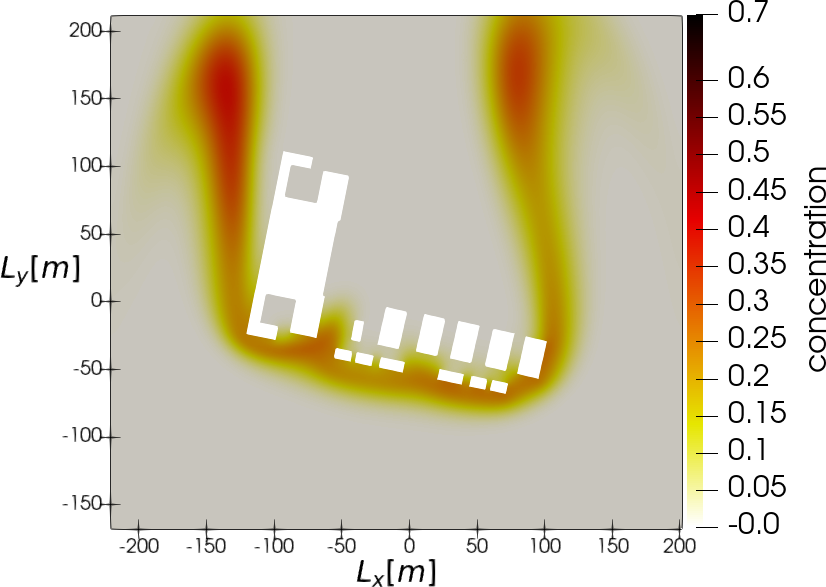}
\end{subfigure}
\caption{Prediction of the contaminant dispersion at ($t=\SI{20}{\s}$) for the two considered sensor configurations (dense left, sparse right).}
\label{fig:eos_pred}
\end{figure}

To illustrate the capabilities and limitations of the proposed method, we consider the campus of the University of the Bundeswehr Munich as a real-world test case. The corresponding wind field is depicted in \autoref{fig:domain}. The computational grid is generated by an automated pipeline that imports building footprints as obstacles directly from OpenStreetMap (OSM) and constructs locally refined triangular meshes to ensure reliable numerical solutions of \autoref{eq:forward_equation} and \autoref{eq:adjoint_equation}, see~\cite{Bonari.2024}. Both partial differential equations are discretized using stabilized linear Lagrange finite elements and implemented within the software framework \fenics~\cite{Baratta.2023}.

As initial parameters, we employ \tcr{the shape function defined in \autoref{eq:ansatz}} with radius $r = \SI{25}{m}$ and threshold $\epsilon = 0.001$. The source trajectory is described by a curve $\gamma : [\SI{0.5}{\s}, \SI{5}{\s}] \rightarrow \R^+ \times \Omega$, illustrated in~\autoref{fig:source}. The forward simulation of the transient source is shown in~\autoref{fig:source}. To demonstrate the performance of the approach, two sensor configurations are considered. The first configuration consists of a dense sensor grid with $149$ sensors measuring concentration values over the observation interval $T = [\SI{0}{\s}, \SI{10}{\s}]$. The measurements are sampled at a rate of $\SI{10}{\Hz}$.
The second configuration uses only $10$ sensors, as depicted in \autoref{fig:eom}. Synthetic measurement data $\measurement$ are generated by adding white Gaussian noise with a signal-to-noise ratio of $\mathrm{SNR} \approx 33.3$. \autoref{fig:eom} additionally illustrates the contaminant distribution at the final measurement time $t=\SI{10}{\s}$. From this instant onward, the prediction of the contaminant concentration is performed.
\begin{figure}
\centering
\includegraphics[width=0.70\linewidth]{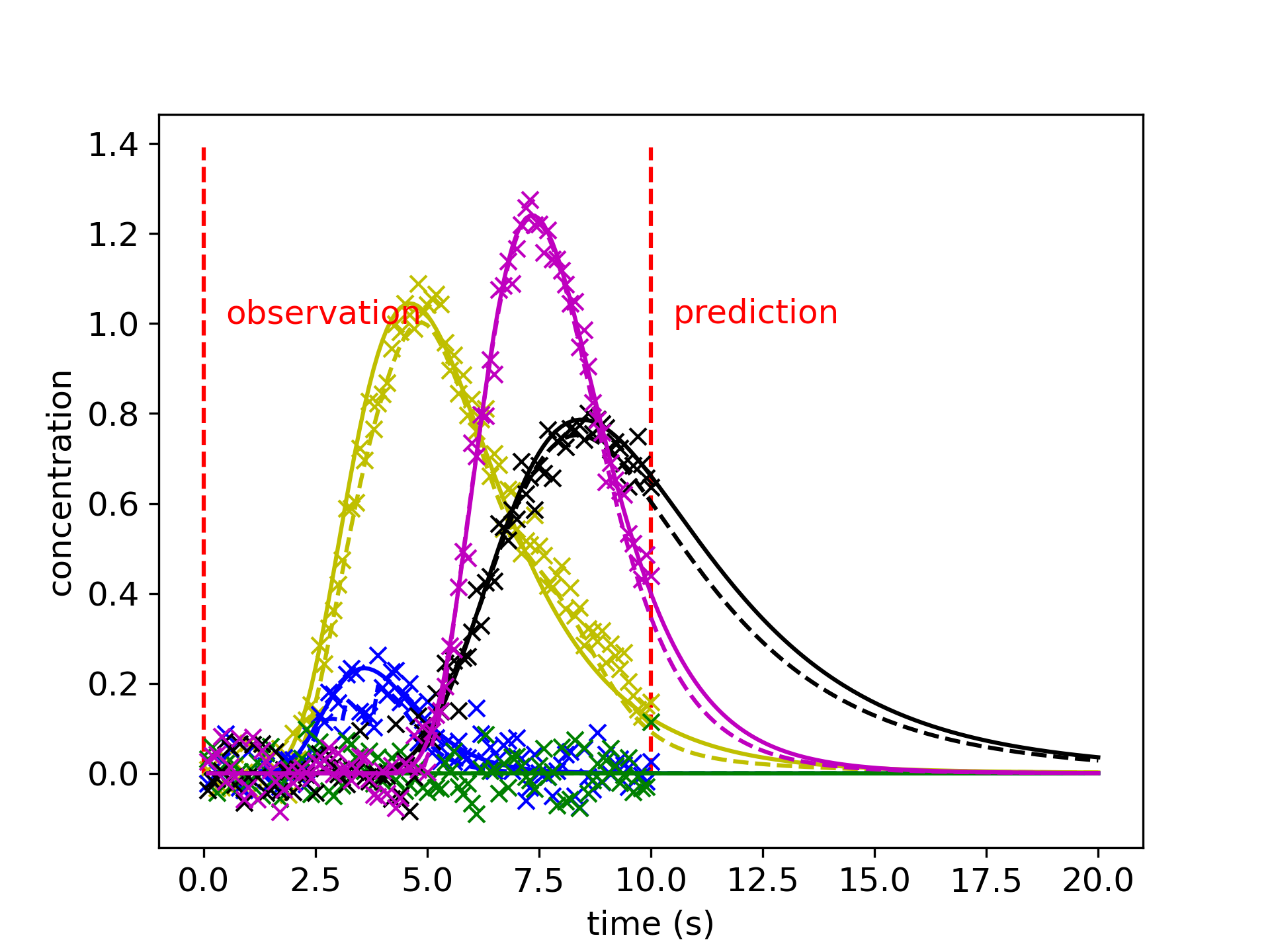}
\caption{Noisy measurements at five sensors (marked by crosses), together with the exact (simulated) concentration (solid line) and the reconstructed concentration (dashed line) obtained using \autoref{alg:pdap} for the coarse sensor grid; cf.~\autoref{fig:eom}.}
\label{fig:time_series}
\end{figure}

We now assess the reconstruction results. \autoref{fig:predstart} displays the reconstructed contaminant distributions for both sensor configurations alongside the true solution (cf.~\autoref{fig:source}). As outlined in the introduction, the algorithm is required to recover not only the spatial source location but also its activation time, which in this example is $t = \SI{0.5}{\s}$. For the dense sensor configuration, the reconstruction closely matches the prescribed parameter curve in both space and time. In contrast, when only $10$ sensors are available, the source location and intensity are recovered only approximately, reflecting the reduced information content of the measurements.
The behavior discussed in the introduction is clearly visible in \autoref{fig:pred04}. At time $t = \SI{3.5}{\s}$, the dense sensor configuration correctly reconstructs the parameter in both space and time. In contrast, for the sparse sensor configuration, the inferred source exhibits only a weak contaminant intensity at this time. The concentrations required to adequately reproduce the sensor signals (cf. \autoref{fig:time_series}) are instead reconstructed by the algorithm at later time instances, reflecting the temporal shift induced by the limited observational information and prior knowledge through the mathematical model. Nevertheless, \autoref{fig:time_series} shows that the measured sensor signals are reproduced very accurately, thereby enabling a reliable prediction of future concentration distributions. At the final release time, $t = \SI{5}{\s}$, the source trajectory is clearly resolved using the dense sensor grid, whereas the sparse configuration again provides only a coarse estimate of the source locations. This effect is shown in \autoref{fig:eosource}. Finally, we examine the predicted contaminant distribution at the end of the measurement interval, $t = \SI{10}{\s}$, and at the final simulation time, for instance $t = \SI{20}{\s}$. In both sensor configurations, the predicted concentration fields show good qualitative agreement with the reference solution. This demonstrates that the proposed method enables reliable forecasts for practical applications, even when the available sensor network is significantly reduced.

\section{Bridging Simulation and Experiment: Model Calibration with 3D Wind Tunnel Data}\label{ssec:wind}

\subsection{Wind Tunnel Experiment}
To test the presented algorithm in a more challenging and realistic scenario, the Authors plan to employ the set-up used in recent experimental studies of gas propagation \cite{hinsen:2024,ruiz:2024}. Here, a measurement campaign has been conducted in a low-speed wind tunnel containing a small-scale model of a set of buildings that represents a hypothetical industrial facility. The cross section of the wind tunnel is subjected to a constant inlet wind velocity of $1.2~\si{\metre\per\second}$, and a synthetic, continuous gas source is placed in the upstream region; theatrical fog with a high concentration of propylene glycol to guarantee thickness and compactness of the plume has been employed. To collect experimental data, a package composed of an array of four sensors travels the whole domain sequentially, resulting in a regular Cartesian acquisition grid of about one cubic decimeter in volume for each cell element. The whole experimental set-up is shown in \autoref{fig:wind_tunnel}.

\subsection{High Fidelity Numerical Simulation Set-up}
The geometrical features of the computational domain are determined by the shape and arrangement of the scaled buildings and the size of the buffer zone to be considered. While the former piece of information is provided in CAD format, the latter represents a set of free parameters that can be categorized in cross-section of the wind tunnel, whose dimensions have been inferred from~\cite{shutin:2024}, and upstream and downstream dimensions, chosen, respectively, as $\lu$ and $\ld$, being $\hmax$ the maximum height of the buildings. According to~\cite{blocken:2015}, these dimensions guarantee a limited effect of the boundaries on the upstream and downstream components of the wind flow.

\begin{figure}
    \centering
    \begin{subfigure}{0.45\textwidth}
        \centering
        \begin{minipage}[c][\textwidth][c]{\textwidth}
            \centering
            \includegraphics[width=\textwidth]{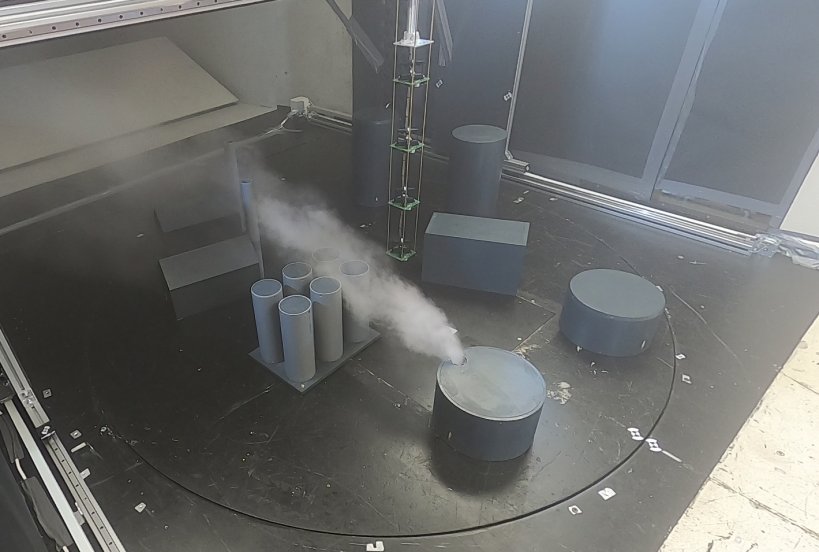}
        \end{minipage}
        \subcaption{}
        \label{subfig:wind_tunnel_a}
    \end{subfigure}
    \hfill
    \begin{subfigure}{0.45\textwidth}
        \centering
        \begin{minipage}[c][\textwidth][c]{\textwidth}
            \centering
            \includegraphics[width=\textwidth]{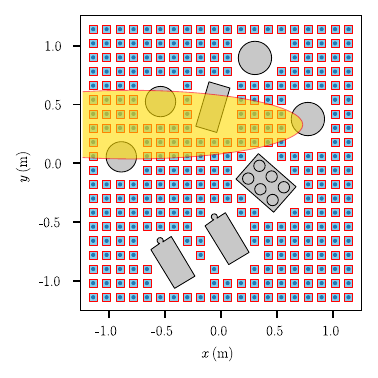}
        \end{minipage}
        \subcaption{}
        \label{subfig:wind_tunnel_b}
    \end{subfigure}
    \caption{Experimental set-up featuring the arrangement of scaled buildings, the synthetic gas plume, and the package of sensors employed to gather data \cite{dlr216980}. Figure~\subref{subfig:wind_tunnel_a} is reproduced with authorization of \cite{hinsen:2024}. Plan view of the sampling process showing all the points acquired at a given elevation of $\overline{z}=0.265\,\si{\meter}$ in Figure~\subref{subfig:wind_tunnel_b}. The difference in domain orientation between Figure~\subref{subfig:wind_tunnel_a} and Figure~\subref{subfig:wind_tunnel_b} are due to data acquired for different wind directions, each obtained rotating the platform underlying the buildings. Data employed in the analysis are related to a domain orientation in accordance with Figure~\subref{subfig:wind_tunnel_b}.}
    \label{fig:wind_tunnel}
\end{figure}

Moreover, a uniform velocity of $1.2\,\si{\metre\per\second}$ is also prescribed at the inlet of the numerical model, no slip conditions are enforced on the surfaces of the buildings, of the ground, and of the wind tunnel. A standard pressure type external-flow condition is prescribed at the outlet. The injection of the contaminant substance is modeled with an additional small inlet surface in correspondence of one of the buildings, cf.~\autoref{fig:wind_tunnel}. Here, a unitary concentration field is prescribed as Dirichlet boundary condition, together with an inlet velocity of $4.0\,\si{\metre\per\second}$. Concentration values are set to zero on all the other surfaces, except for the global outlet, where a zero Neumann boundary condition for the concentration field is enforced. The mesh is built using the Ansys-Fluent mosaic Poly-Hexcore watertight meshing workflow, a hybrid meshing scheme that combines \emph{hexahedral} and \emph{polyhedral} elements to allow an optimal structured mesh in the core region of the geometry, while relying on the adaptability of polyhedral elements where regular meshing is not possible.

A staggered approach is employed to solve the dispersion problem. First, the wind field alone gets evaluated, based on the aforementioned conditions and employing a~\emph{k-$\omega$} set of equations to model turbulence. In a subsequent solution step, the AD equation alone is solved, providing the wind field as a known parameter and resulting in a steady state contaminant cloud that aims at approximating the dispersion conditions of the experimental dispersion process. This modular approach has the advantage that each component can be modified independently, e.g., the contaminant transport problem can be easily extended towards transient conditions.

\subsection{Identifying Model Parameters Best Explaining the Experimental Data}

\begin{figure}
    \centering
    \includegraphics[width=.4\textwidth]{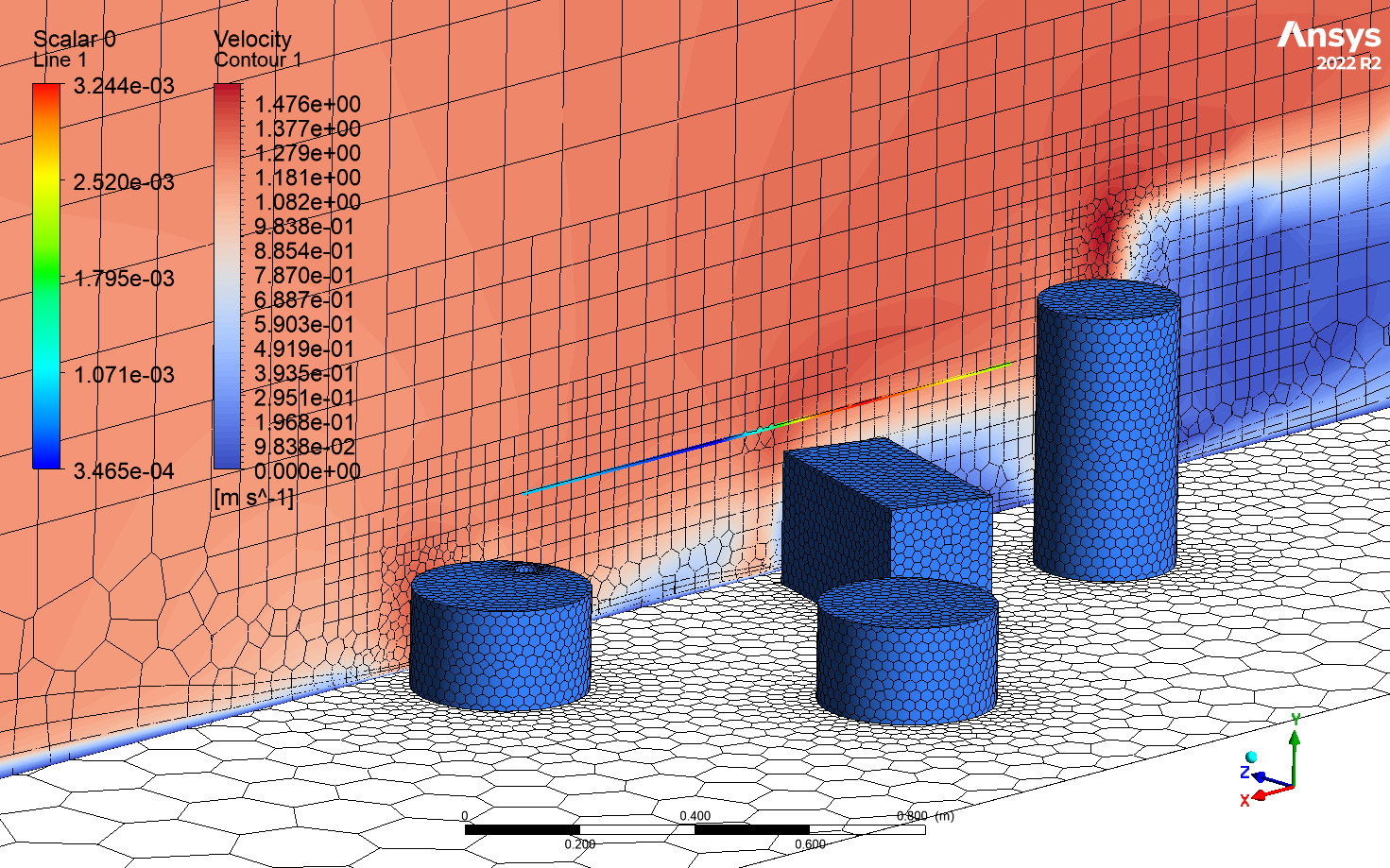}    
    \quad
    \includegraphics[width=.4\textwidth]{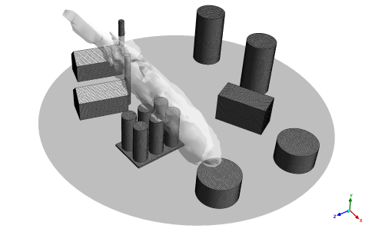}
    \caption{Excerpt of the computational domain featuring a portion of the buildings arrangement, contour values of wind on an exemplary section cut and the line where the concentration values have been sampled to compare with results from the experiment. Simulated concentration values are superposed to the line element as additional plot (left). Simulated contaminant dispersion in qualitative agreement with wind tunnel experiment (~\autoref{fig:wind_tunnel}) (right).}
    \label{fig:fff_line}
\end{figure}

\begin{figure}
    \centering
    \includegraphics[width=\textwidth]{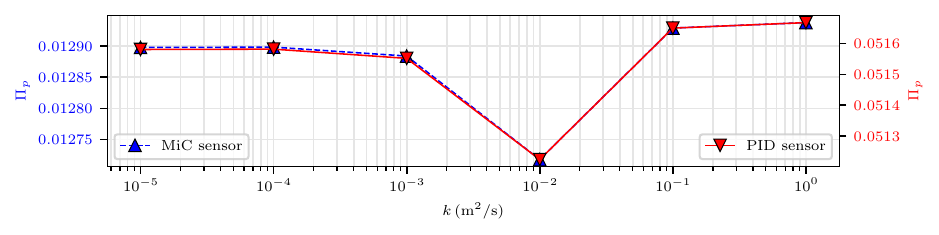}
    \caption{Numerical evaluation of the objective expressions $\bPi$ for realizations of $k = \kappa_i$. The subscripts~\emph{(m)} and~\emph{(p)} refer to the two experimental sets under examination, the first related to experimental data acquired through a SGX Sensortech~\emph{MiCS-5524} MOX sensor and the second through an Alphasense~\emph{PID-AH2} PID sensor.}
    \label{fig:objective}
\end{figure}

In the calibration stage, a tentative is made to acquire a realistic\footnote{Here, the term \emph{realistic} has the meaning of \emph{best explaining the experimental data given the AD model assumptions.}} value of the diffusion coefficient $\kappa$. To this purpose, a parameter space is defined considering different values spanning several orders of magnitude, collected in the array:
\begin{equation}
    \vec{\kappa}=\{10^i,\,i=-5,\dots,0\}\in\mathbb{R}^6,
\end{equation}

and a model run is performed for each of the $\kappa_i$ values. With reference to the numerical solutions, concentration values are virtually sampled on several points of the computational domain, coincident with a subset of sensor readings in the experimental setting. More specifically, the acquisition is performed on $\np=11$ points equally spaced on a line parallel to the global $x$ direction, with starting point $\bx=[0.54\,\si{\meter},0.18\,\si{\meter},0.27\,\si{\meter}]$ and length $l=1.2\,\si{\meter}$, cf.~\autoref{fig:fff_line}.
The values obtained at these points have been compared to the two sets of concentration readings resulting from the experiments, each of the two sets being related to a different type of sensor employed \cite{shutin:2024}. The value to be used is then found, among the tentative ones, as the minimizer of the discrete cost objective expression:
\begin{equation}
    \Pi_{i,j} = \frac{1}{N}\sum_{l=1}^{\np}\Bigl(\cs_l(\kappa_i)-\ce_{l,j}\Bigl)^2,\quad j={m,p},
\end{equation}
where the subscript $j$ identifies the two types of sensors employed in the experiment, $l$ the collocation point in the domain, and the arrays $\bcs$ and $\bce_{(1,2)}$ collect the numerically computed (simulated) concentration values and measured concentration readings in the wind tunnel experiment, respectively. The values of the objective expressions are collected in the arrays $\bPim$ and $\bPip$, respectively, and represented in \autoref{fig:objective}. This results show that a value of $\kappa \approx10^{-2}\,\si{\metre^2\per\second}$ is the best fit to approximate the data for both sets of sensor readings available. The result of the preliminary study shows a good agreement between the identified optimal value of the diffusion coefficient and the one identified in~\cite{hinsen:2024}.

\section{A Strategy for real-time data coupling and digital twin integration}\label{sec:dtwin}
So far, a simple parameter sweep was performed to determine the diffusion coefficient $\kappa$  that represents best the observed situation in the wind tunnel. In that initial study all observations were assumed to be available simultaneously for model calibration. In a CBRN incident, however, data might become available only sequentially, which precludes batch calibration and necessitates an online‐identification framework.

To address this, we plan to employ sequential Bayesian inference (SBI) methods that update the parameter estimates as new measurements arrive, as for instance Kalman filters. For parameter estimation based on SBI there are two common strategies: the augmented-state approach and the dual-filtering approach. The augmented‑state approach expands the state vector to include the unknown parameters and treats them as additional dynamical variables \cite{Asch.2022}. Alternatively, a dual‑filtering scheme can be applied, where one filter propagates the physical state while a second filter estimates the parameters \cite{Chebbi.2025}.  
Because the advection–diffusion model yields a high‑dimensional state vector \cite{10773899} and its numerical evaluation is computationally expensive, we propose to use Ensemble Kalman Filters (EnKF) in combination with one of the parameter estimation strategies described above. EnKFs are particularly suited to high‑dimensional, potentially nonlinear problems, as unlike the classical Kalman filter they avoid repeated inversion of large covariance matrices and do not require explicit linearization of the model \cite{Asch.2022}. Their effectiveness in comparable settings has already been demonstrated by \cite{10773899}. Comprehensive introductions and in-depth treatments of data assimilation can be found in \cite{Asch.2022, Asch.2016}. 
In addition, rapid parameter estimation in such high‑dimensional contexts demands not only an efficient data‑assimilation scheme but also fast evaluations of the forward model. This can be accomplished, for instance, by leveraging surrogate techniques from physics‑informed machine learning \cite{Griese.2025} or by employing reduced‑order modeling \cite{Asch.2022}.

Moreover, the wind field was assumed to be constant during the source identification task. If the wind conditions change during the course of the incident, a complete reevaluation of the high-fidelity (HF) CFD model is computationally too expensive. To overcome this limitation, a parameterized reduced-order model (ROM) can be derived from the corresponding HF model as an accurate approximation that can be evaluated in only a fraction of the time required by the full-order simulation. Model order reduction methods are often distinguished into intrusive and non-intrusive approaches, depending on whether access to the high-fidelity operators is required. While intrusive methods are purely physics-based, they require access to the HF operators, which can be problematic when closed-source or legacy solvers are used. Examples of such methods can be found in \cite{Quarteroni.2015,Rozza.2022,Hesthaven.2026}. Non-intrusive methods are typically purely data-driven, which comes with its own drawbacks; however, they often provide substantial speed-ups and enable applications in cases where no access to the high-fidelity operators is available. Such methods can be found in \cite{Yu.2019,Fresca.2022,Vinuesa.2022}. The prediction of wind fields under varying ambient wind conditions for the considered application is treated, e.g., in \cite{kühn2026} and references therein.

To have pre-computed data readily available in case of an incident, we plan to connect our simulations with an automatically generated urban digital twin \cite{10732278}. Extending a component-based digital twin framework (\cite{10407185}) towards database-centric simulation data management, we prepare the proposed algorithm for interactive digital twin applications in disaster management. 

\section{Summary and Outlook}
\begin{figure}
\centering
\includegraphics[width=.8\textwidth]{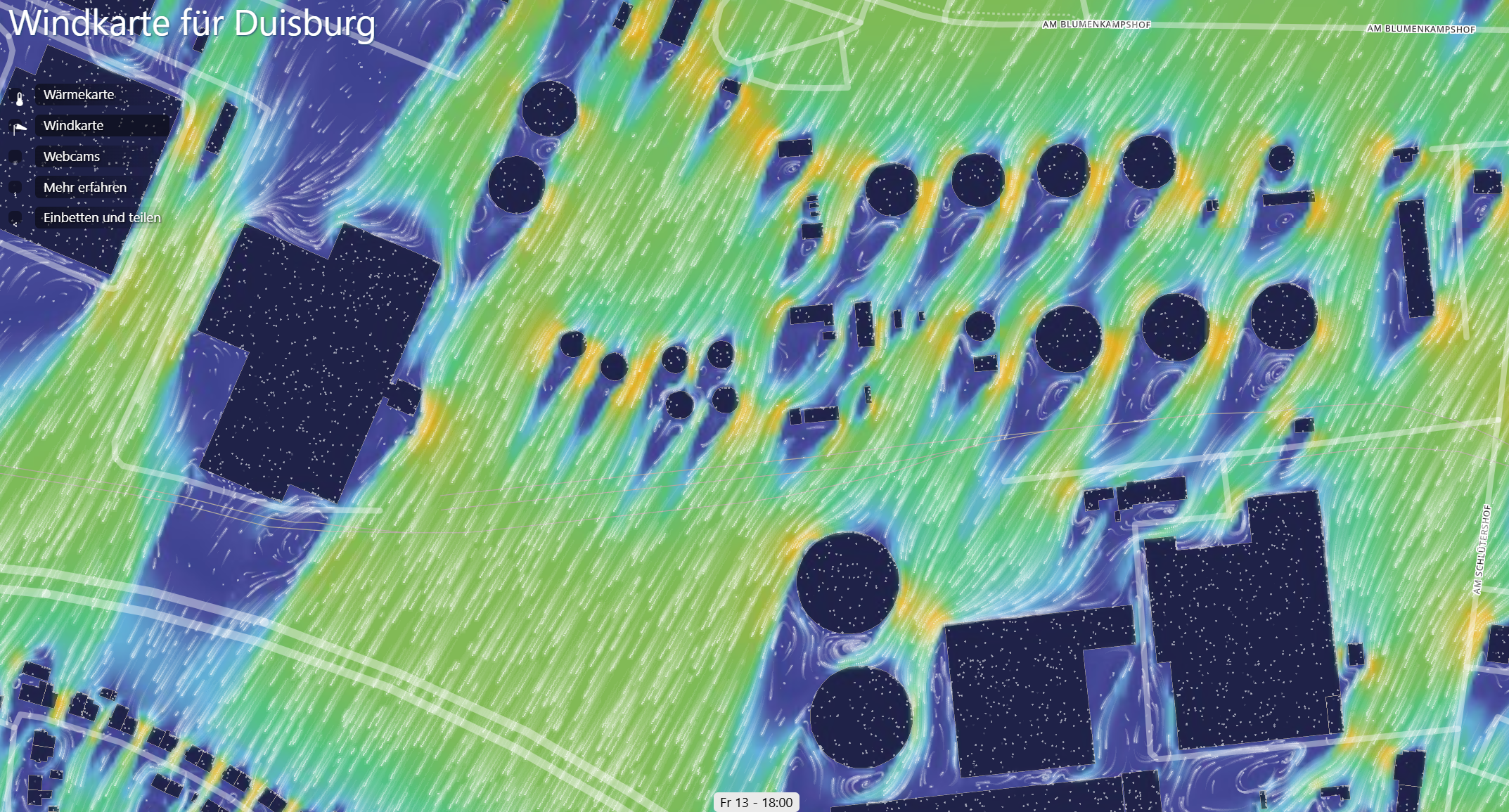}
\caption{Wind field in the area of the port of Duisburg~\cite{WindDuisburg} as example input for the contaminant source identification algorithm.}
\label{fig:windfig}
\end{figure}
This work presents a novel algorithm for the rapid identification of moving contaminant sources from sparse sensor measurements governed by an advection–diffusion model. 
Owing to the high efficiency of the proposed method (see~\cite{MATTUSCHKA2026118854}), fewer than $100$ iterations of \autoref{alg:pdap} are required to accurately reconstruct the sources. Each iteration involves the solution of the forward problem~\autoref{eq:forward_equation} and the corresponding adjoint problem~\autoref{eq:adjoint_equation}. In contrast, sampling-based approaches such as the Metropolis–Hastings algorithm typically require on the order of $10{.}000$ or more PDE solves, making highly efficient and intelligent surrogate modeling indispensable (see, e.g.,~\cite{AMMAR2026121873}).

The computational efficiency of the present formulation enables, for the first time, the identification of transient sources on large spatial scales. In the numerical experiments, convergence is achieved after only seventeen iterations for a dense sensor grid and $36$ iterations for a coarse grid. Although each iteration of \autoref{alg:pdap} requires $N_T$ forward solves~\autoref{eq:forward_equation} and one adjoint solve~\autoref{eq:adjoint_equation}, the forward computations are fully independent and therefore straightforward to parallelize. Consequently, the overall computational cost is comparable to that of identifying an initial condition. Moreover, as demonstrated in~\cite{MATTUSCHKA2026118854}, the method remains effective for multiple simultaneous sources while requiring only a limited number of PDE solves.

Nevertheless, the numerical examples also reveal a limitation of the current approach: accurate reconstruction of source activation times deteriorates when the sensor density is reduced. Moreover, in this setting it is not possible to derive a reasonable optimal sensor placement (optimal design of experiments), cf.~\cite{Huynh_2024,MattuschkaGoal}. To address this issue, additional methodological extensions, analogous to those introduced for the heat conduction equation, would be required, cf.~\cite{4f98d6b8773f48d5ba2099503bed735c, Gong.2025}

To enable the application of the proposed method in realistic scenarios, the previous \autoref{ssec:wind} outlined the necessary methodological and modeling extensions. The next objective is to validate the approach in the context of real-world wind fields, for example using an urban setting such as the city of Duisburg (see \autoref{fig:windfig}), and to support the numerical results with experimental data. Possible strategies for the required model calibration were discussed in previous \autoref{sec:dtwin} and will form a central component of future research. In particular, the systematic validation and adjustment of the flow and transport models are essential steps toward transferring the method from a controlled computational environment to operational real-world applications. In a further application scenario an attack with multiple drones could be considered.
Assuming that drones are detected image-based~\cite{lenhard2025syndronevision}, the algorithm could also help to identify critical drones that are in fact contaminant sources and distinguish them from others acting as camouflage.

\section{Acknowledgements}
We thank Lisa Kühn and Philip Franz for providing insights on model-order reduction and filtering techniques. AP gratefully acknowledges the funding by dtec.bw - Digitalization and Technology Research Center of the Bundeswehr (project RISK.twin). dtec.bw is funded by the European Union - NextGenerationEU.

\printbibliography

@article{shutin:2024,
	title = {Gas source localization from real-world spatial in-situ concentration and wind measurements},
	author = {Shutin, D and Munoz, C and Wiedemann, T and Hinsen, P and Ruiz, V Prieto and Zhang, S and Lilienthal, A J and Fan, H},
	langid = {english},
	year = {2025},
	doi= {10.21227/x0sf-ad36} 
}

@article{blocken:2015,
	title = {Computational Fluid Dynamics for urban physics: Importance, scales, possibilities, limitations and ten tips and tricks towards accurate and reliable simulations},
	volume = {91},
	issn = {03601323},
	doi = {10.1016/j.buildenv.2015.02.015},
	shorttitle = {Computational Fluid Dynamics for urban physics},
	pages = {219--245},
	journaltitle = {Building and Environment},
	shortjournal = {Building and Environment},
	author = {Blocken, Bert},
	urldate = {2024-02-14},
	date = {2015-09},
	langid = {english},
}

@inproceedings{hinsen:2024,
	title = {Experimental Study of Gas Propagation: Parameter Identification and Analysis in a Wind Tunnel},
	rights = {https://doi.org/10.15223/policy-029},
	doi = {10.1109/isoen61239.2024.10556306},
	shorttitle = {Experimental Study of Gas Propagation},
    eventtitle = {2024 {IEEE} International Symposium on Olfaction and Electronic Nose ({ISOEN})},
	pages = {1--3},
	booktitle = {2024 {IEEE} International Symposium on Olfaction and Electronic Nose ({ISOEN})},
	publisher = {{IEEE}},
	author = {Hinsen, Patrick and Wiedemann, Thomas and Ruiz, Victor Scott Prieto and Shutin, Dmitriy and Lilienthal, Achim J.},
	urldate = {2025-07-08},
	date = {2024-05-12},
	langid = {english},
	}

@inproceedings{ruiz:2024,
	location = {Grapevine, {TX}, {USA}},
	title = {Gas Source Localization Using Physics-Guided Neural Networks},
	rights = {https://doi.org/10.15223/policy-029},
	doi = {10.1109/isoen61239.2024.10556061},
    eventtitle = {2024 {IEEE} International Symposium on Olfaction and Electronic Nose ({ISOEN})},
	pages = {1--3},
	booktitle = {2024 {IEEE} International Symposium on Olfaction and Electronic Nose ({ISOEN})},
	publisher = {{IEEE}},
	author = {Ruiz, Victor Prieto and Hinsen, Patrick and Wiedemann, Thomas and Shutin, Dmitriy and Christof, Constantin},
	urldate = {2025-07-08},
	date = {2024-05-12},
	langid = {english},
}

@article{.,
 title = {paper{\_}automated{\_}model{\_}generation{\_}validated}
}

@book{Asch.2016,
 author = {Asch, Marc and Asch, Mark and Bocquet, Marc and Nodet, Ma{\"e}lle},
 year = {2016},
 title = {Data assimilation: Methods, algorithms, and applications},
 price = {Paperback : circa USD 84.00},
 address = {Philadelphia},
 volume = {11},
 publisher = {{siam Society for Industrial and Applied Mathematics}},
 isbn = {9781611974539},
 series = {Fundamentals of algorithms},
 doi = {10.1137/1.9781611974546}
}

@article{Griese.2025,
 author = {Griese, Franziska and Hoppe, Fabian and R{\"u}ttgers, Alexander and Knechtges, Philipp},
 year = {2025},
 title = {Preconditioned FEM-based neural networks for solving incompressible fluid flows and related inverse problems},
 pages = {116663},
 volume = {469},
 journal = {Journal of Computational and Applied Mathematics},
 doi = {10.1016/j.cam.2025.116663}
}

@Article{Chebbi.2025,
AUTHOR = {Chebbi, Amal and Franchek, Matthew A. and Grigoriadis, Karolos},
TITLE = {Simultaneous State and Parameter Estimation Methods Based on Kalman Filters and Luenberger Observers: A Tutorial \& Review},
JOURNAL = {Sensors},
VOLUME = {25},
YEAR = {2025},
NUMBER = {22},
ARTICLE-NUMBER = {7043},
PubMedID = {41305249},
ISSN = {1424-8220},
DOI = {10.3390/s25227043}
}

@book{Asch.2022,
 author = {Asch, Mark},
 year = {2022},
 title = {A toolbox for digital twins: From model-based to data-driven},
 address = {Philadelphia},
 volume = {6},
 publisher = {{Society for Industrial and Applied Mathematics}},
 isbn = {9781611976977},
 series = {Mathematics in industry},
 doi = {10.1137/1.9781611976977}
}

@article{Biccari.2023,
 author = {Biccari, Umberto and Song, Yongcun and Yuan, Xiaoming and Zuazua, Enrique},
 year = {2023},
 title = {A two-stage numerical approach for the sparse initial source identification of a diffusion--advection equation},
 pages = {095003},
 volume = {39},
 number = {9},
 journal = {Inverse Problems},
 doi = {10.1088/1361-6420/ace548}
}

@inproceedings{Bonari.2024,
 author = {Bonari, Jacopo and K{\"u}hn, Lisa and von Danwitz, Max and Popp, Alexander},
 title = {Towards Real-Time Urban Physics Simulations with Digital Twins},
 pages = {18--25},
 publisher = {IEEE},
 isbn = {979-8-3315-2721-1},
 booktitle = {2024 28th International Symposium on Distributed Simulation and Real Time Applications (DS-RT)},
 year = {2024},
 doi = {10.1109/DS-RT62209.2024.00013}
}

@article{Boris.2002,
 author = {Boris, J.},
 year = {2002},
 title = {The threat of chemical and biological terrorism: preparing a response},
 pages = {22--32},
 volume = {4},
 number = {2},
 issn = {15219615},
 journal = {Computing in Science {\&} Engineering},
 doi = {10.1109/5992.988644}
}

@article{Brooks.1982,
 author = {Brooks, Alexander N. and Hughes, Thomas J.R.},
 year = {1982},
 title = {Streamline upwind/Petrov-Galerkin formulations for convection dominated flows with particular emphasis on the incompressible Navier-Stokes equations},
 pages = {199--259},
 volume = {32},
 number = {1-3},
 journal = {Computer Methods in Applied Mechanics and Engineering},
 doi = {10.1016/0045-7825(82)90071-8}
}

@article{Danwitz.2023,
 author = {von Danwitz, Max and Voulis, Igor and Hosters, Norbert and Behr, Marek},
 year = {2023},
 title = {Time--continuous and time--discontinuous space--time finite elements for advection--diffusion problems},
 pages = {3117--3144},
 volume = {124},
 number = {14},
 journal = {International Journal for Numerical Methods in Engineering},
 doi = {10.1002/nme.7241}
}

@article{Milzarek.2014,
 author = {Milzarek, Andre and Ulbrich, Michael},
 year = {2014},
 title = {A Semismooth Newton Method with Multidimensional Filter Globalization for $l_1$-Optimization},
 pages = {298--333},
 volume = {24},
 number = {1},
 issn = {1052-6234},
 journal = {SIAM Journal on Optimization},
 doi = {10.1137/120892167}
}

@article{Monge.2020,
 author = {Monge, Azahar and Zuazua, Enrique},
 year = {2020},
 title = {Sparse source identification of linear diffusion--advection equations by adjoint methods},
 pages = {104801},
 volume = {145},
 journal = {Systems {\&} Control Letters},
 doi = {10.1016/j.sysconle.2020.104801}
}

@incollection{Patnaik.2012,
 author = {Patnaik, Gopal and Boris, Jay P. and Grinstein, Fernando F. and Iselin, John P. and Hertwig, Denise},
 title = {Large Scale Urban Simulations with FCT},
 pages = {91--117},
 publisher = {{Springer Netherlands}},
 isbn = {978-94-007-4037-2 },
 series = {Scientific Computation},
 editor = {Kuzmin, Dmitri and L{\"o}hner, Rainald and Turek, Stefan},
 booktitle = {Flux-Corrected Transport},
 year = {2012},
 address = {Dordrecht},
 doi = {10.1007/978-94-007-4038-9_4 }
}

@article{Pieper.2021,
 author = {Pieper, Konstantin and Walter, Daniel},
 year = {2021},
 title = {Linear convergence of accelerated conditional gradient algorithms in spaces of measures},
 pages = {38},
 volume = {27},
 issn = {1292-8119},
 journal = {ESAIM: Control, Optimisation and Calculus of Variations},
 doi = {10.1051/cocv/2021042 }
}

@article{Villa.2021,
 author = {Villa, Umberto and Petra, Noemi and Ghattas, Omar},
 year = {2021},
 title = {hIPPYlib},
 pages = {1--34},
 volume = {47},
 number = {2},
 journal = {ACM Transactions on Mathematical Software},
 doi = {10.1145/3428447}
}

@article{Vinuesa.2022,
 author = {Vinuesa, Ricardo and Brunton, Steven L.},
 year = {2022},
 title = {Enhancing computational fluid dynamics with machine learning},
 pages = {358--366},
 volume = {2},
 number = {6},
 journal = {Nature computational science},
 doi = {10.1038/s43588-022-00264-7}
}

@article{Yu.2019,
 author = {Yu, Jian and Yan, Chao and Guo, Mengwu},
 year = {2019},
 title = {Non-intrusive reduced-order modeling for fluid problems: A brief review},
 pages = {5896--5912},
 volume = {233},
 number = {16},
 journal = {Proceedings of the Institution of Mechanical Engineers, Part G: Journal of Aerospace Engineering},
 doi = {10.1177/0954410019890721}
}

@article{Casas.2019,
 author = {Casas, E. and Kunisch, K.},
 year = {2019},
 title = {Using sparse control methods to identify sources in linear diffusion-convection equations},
 pages = {114002},
 volume = {35},
 number = {11},
 issn = {0266-5611},
 journal = {Inverse Problems},
 doi = {10.1088/1361-6420/ab331c}
}

@article{Leykekhman.2020,
 author = {Leykekhman, Dmitriy and Vexler, Boris and Walter, Daniel},
 year = {2020},
 title = {Numerical analysis of sparse initial data identification for parabolic problems},
 pages = {1139--1180},
 volume = {54},
 number = {4},
 issn = {0764-583X},
 journal = {ESAIM: Mathematical Modelling and Numerical Analysis},
 doi = {10.1051/m2an/2019083}
}

@article{Elman.2020,
 author = {Elman, Howard C. and Su, Tengfei},
 year = {2020},
 title = {A low-rank solver for the stochastic unsteady Navier--Stokes problem},
 pages = {112948},
 volume = {364},
 issn = {0045-7825},
 journal = {Computer Methods in Applied Mechanics and Engineering},
 doi = {10.1016/j.cma.2020.112948}
}

@article{MATTUSCHKA2026118854,
title = {Sparse source identification in transient advection-diffusion problems with a primal-dual-active-point strategy},
journal = {Computer Methods in Applied Mechanics and Engineering},
volume = {454},
pages = {118854},
year = {2026},
issn = {0045-7825},
doi = {10.1016/j.cma.2026.118854},
author = {Marco Mattuschka and Daniel Walter and Max {von Danwitz} and Alexander Popp}
}

@article{MattuschkaGoal, title={Goal-oriented optimal sensor placement for PDE-constrained inverse problems in crisis management}, DOI={10.24423/cames.2026.1887}, abstractNote={
This paper presents a novel framework for goal-oriented optimal static sensor placement and dynamic sensor steering in PDE-constrained inverse problems, utilizing a Bayesian approach accelerated by low-rank approximations. The framework is applied to airborne contaminant tracking, extending recent dynamic sensor steering methods to complex geometries for computational efficiency. A C-optimal design criterion is employed to strategically place sensors, minimizing uncertainty in predictions. Numerical experiments validate the approach’s effectiveness for source identification and monitoring, highlighting its potential for real-time decision-making in crisis management scenarios.
}, journal={Computer Assisted Methods in Engineering and Science}, author={Mattuschka, Marco and An der Lan, Noah and von Danwitz, Max and Wolff, Daniel and Alexander , Popp}, year={2026}, month={2} }

@incollection{Danwitz.2024,
 author = {{von Danwitz}, Max and Bonari, Jacopo and Franz, Philip and K{\"u}hn, Lisa and Mattuschka, Marco and Popp, Alexander},
 title = {Contaminant Dispersion Simulation in a Digital Twin Framework for Critical Infrastructure Protection},
 pages = {1--12},
 year={2024},
 booktitle = {Towards Digital Twins for Infrastructures},
 doi = {10.23967/eccomas.2024.301},
}

@INPROCEEDINGS{10773899,
  author={Gioia, Daniele Giovanni and Bonari, Jacopo and Lichte, Daniel and Popp, Alexander},
  booktitle={2024 Sensor Data Fusion: Trends, Solutions, Applications (SDF)}, 
  title={Sequential Drone Routing for Data Assimilation on a 2D Airborne Contaminant Dispersion Problem}, 
  year={2024},
  volume={},
  number={},
  pages={1-8},
  doi={10.1109/SDF63218.2024.10773899}}

@article{4f98d6b8773f48d5ba2099503bed735c,
title = "Measure valued directional sparsity for parabolic optimal control problems",
author = "Karl Kunisch and Konstantin Pieper and Boris Vexler",
note = "Publisher Copyright: {\textcopyright} 2014 Society for Industrial and Applied Mathematics",
year = "2014",
doi = "10.1137/140959055",
language = "English",
volume = "52",
pages = "3078--3108",
journal = "SIAM Journal on Control and Optimization",
issn = "0363-0129",
publisher = "Society for Industrial and Applied Mathematics Publications",
number = "5",
}

@article{Gong.2025,
 author = {Gong, Wei and Liang, Dongdong},
 year = {2025},
 title = {Analysis and approximation to parabolic optimal control problems with measure-valued controls in time},
 pages = {2},
 volume = {31},
 issn = {1292-8119},
 journal = {ESAIM: Control, Optimisation and Calculus of Variations},
 doi = {10.1051/cocv/2024085}
}

@book{Quarteroni.2015,
 author = {Quarteroni, Alfio and Manzoni, Andrea and Negri, Federico},
 year = {2015},
 title = {Reduced basis methods for partial differential equations: An introduction},
 address = {Cham},
 volume = {92},
 publisher = {Springer},
 isbn = {978-3-319-15431-2},
 series = {UNITEXT La Matematica per il 3+2},
 doi = {10.1007/978-3-319-15431-2}
}

@book{Rozza.2022,
 year = {2022},
 title = {Advanced Reduced Order Methods and Applications in Computational Fluid Dynamics},
 address = {Philadelphia, PA},
 publisher = {{Society for Industrial and Applied Mathematics}},
 isbn = {978-1-61197-724-0},
 editor = {Rozza, Gianluigi and Stabile, Giovanni and Ballarin, Francesco},
 doi = {10.1137/1.9781611977257}
}

@misc{Hesthaven.2026,
 author = {Hesthaven, Jan S. and Peherstorfer, Benjamin and Unger, Benjamin},
 date = {2026},
 title = {Nonlinear model reduction for transport-dominated problems},
 publisher = {arXiv},
 doi = {10.48550/arXiv.2602.01397}
}

@misc{Baratta.2023,
 author = {Baratta, Igor A. and Dean, Joseph P. and Dokken, J{\o}rgen S. and Habera, Michal and Hale, Jack S. and Richardson, Chris N. and Rognes, Marie E. and Scroggs, Matthew W. and Sime, Nathan and Wells, Garth N.},
 year = {2023},
 title = {{DOLFINx}: The next generation {FEniCS} problem solving environment},
 publisher = {Zenodo},
 doi = {10.5281/zenodo.10447666}
}

@article{Fresca.2022,
 author = {Fresca, Stefania and Manzoni, Andrea},
 year = {2022},
 title = {POD-DL-ROM: Enhancing deep learning-based reduced order models for nonlinear parametrized PDEs by proper orthogonal decomposition},
 pages = {114181},
 volume = {388},
 issn = {00457825},
 journal = {Computer Methods in Applied Mechanics and Engineering},
 doi = {10.1016/j.cma.2021.114181}
}

@article{WANG2025180041,
title = {Deep learning based aerosol particle classification for the detection of ship emissions},
journal = {Science of The Total Environment},
volume = {994},
pages = {180041},
year = {2025},
issn = {0048-9697},
doi = {10.1016/j.scitotenv.2025.180041},
author = {Guanzhong Wang and Heinrich Ruser and Julian Schade and Seongho Jeong and Johannes Passig and Ralf Zimmermann and Günther Dollinger and Thomas Adam}
}

@inproceedings{spie:133423fc396d96d6105c143c6a891646b5384334,
	title = "Combining standoff tomography with point detection: a game changer for the identification of airborne toxic chemicals",
	booktitle = "Proc. SPIE",
	author = "Wilsenack, Frank and Allers, Maria and Meyer, Fabian and Wolf, Thomas and Tjärnhage, Torbjörn and Landström, Lars and Ficks, Arne",
	issn = "0277-786X",
	isbn = "9781510674301",
	volume = "13056",
	pages = "130560Z",
	publisher = "SPIE; ",
	doi = {10.1117/12.3013427},
	year = {2024},
}

@inproceedings{10.1117/12.692922,
author = {Roland Harig and J{\"o}rn Gerhard and Ren{\'e} Braun and Chris Dyer and Ben Truscott and Richard Moseley},
title = {{Remote detection of gases and liquids by imaging Fourier transform spectrometry using a focal plane array detector: first results}},
volume = {6378},
booktitle = {Chemical and Biological Sensors for Industrial and Environmental Monitoring II},
editor = {Steven D. Christesen and Arthur J. Sedlacek III and James B. Gillespie and Kenneth J. Ewing},
organization = {International Society for Optics and Photonics},
publisher = {SPIE},
pages = {637816},
year = {2006},
doi = {10.1117/12.692922},
}

@article{MADADELAHI2025117099,
title = {Electrochemical sensors: Types, applications, and the novel impacts of vibration and fluid flow for microfluidic integration},
journal = {Biosensors and Bioelectronics},
volume = {272},
pages = {117099},
year = {2025},
issn = {0956-5663},
doi = {10.1016/j.bios.2024.117099},
author = {Masoud Madadelahi and Fabian O. Romero-Soto and Rudra Kumar and Uriel Bonilla Tlaxcala and Marc J. Madou}
}

@misc{kühn2026,
      title={Intrusive and Non-Intrusive Model Order Reduction for Airborne Contaminant Transport: Comparative Analysis and Uncertainty Quantification}, 
      author={Lisa Kühn and Jacopo Bonari and Max von Danwitz and Alexander Popp},
      year={2026},
      doi={10.48550/arXiv.2602.21996 }, 
}

@article{AMMAR2026121873,
title = {Inverse identification of a passive scalar source of pollution in large areas},
journal = {Atmospheric Environment},
volume = {372},
pages = {121873},
year = {2026},
issn = {1352-2310},
doi = {10.1016/j.atmosenv.2026.121873},
author = {Amine Ammar and Francisco Chinesta}
}

@article{Huynh_2024,
doi = {10.1088/1361-6420/ad2cf8},
year = {2024},
publisher = {IOP Publishing},
volume = {40},
number = {5},
pages = {055007},
author = {Huynh, Phuoc-Truong and Pieper, Konstantin and Walter, Daniel},
title = {Towards optimal sensor placement for inverse problems in spaces of measures},
journal = {Inverse Problems}
}

@article{dlr216980,
       publisher = {Nature Publishing Group},
          author = {Hinsen, Patrick and Wiedemann, Thomas and Shutin, Dmitriy and Lilienthal, Achim},
         journal = {Scientific Data},
            year = {2026},
           title = {High-Resolution Wind Tunnel Dataset of Gas Sensor Responses to Vapor Plumes in Scale Model Landscapes},
            issn = {2052-4463},
             url = {https://elib.dlr.de/216980/}
}

@inproceedings{lenhard2025syndronevision,
  title={SynDroneVision: A synthetic dataset for image-based drone detection},
  author={Lenhard, Tamara R and Weinmann, Andreas and Franke, Kai and Koch, Tobias},
  booktitle={2025 IEEE/CVF Winter Conference on Applications of Computer Vision (WACV)},
  pages={7637--7647},
  year={2025},
  organization={IEEE},
  doi={10.1109/WACV61041.2025.00742}
}

@online{WindDuisburg,
  title     = {Windkarte für Duisburg},
  author    = {meteoblue, A Windy.com Company},
  year      = 2026,
  url       = {https://www.meteoblue.com/de/products/cityclimate/heatmaps/duisburg?map=wind#16.62/51.439359/6.734443},
  urldate   = {2026-02-13}
}

@INPROCEEDINGS{10732278,
  author={Cao, Susanna and Franke, Kai and Koch, Tobias and Lenhard, Tamara},
  booktitle={2024 Annual Modeling and Simulation Conference (ANNSIM)}, 
  title={Automated Generation Of Urban Digital Shadows Using Open Data}, 
  year={2024},
  volume={},
  number={},
  pages={1-13},
  doi={10.23919/ANNSIM61499.2024.10732278}}

@INPROCEEDINGS{10407185,
  author={Franke, Kai and Stürmer, J. Marius and Koch, Tobias},
  booktitle={2023 Winter Simulation Conference (WSC)}, 
  title={Automated Simulation and Virtual Reality Coupling for Interactive Digital Twins}, 
  year={2023},
  volume={},
  number={},
  pages={2615-2626},
  doi={10.1109/WSC60868.2023.10407185}}

\end{document}